\documentclass[11pt]{article}

\usepackage{hyperref}

\usepackage{array} % for the table at the end
\usepackage{amsfonts}
\usepackage{amsmath}
\usepackage{amssymb}
\usepackage{graphicx}
\usepackage{color}
\usepackage{esint}    % for averages

\usepackage{enumitem} % for enumerate/itemize margins

\newtheorem{thm}{Theorem}[section]

\newtheorem{lem}[thm]{Lemma}
\newtheorem{prop}[thm]{Proposition}
\newtheorem{defn}[thm]{Definition}

\numberwithin{equation}{section}

\newcommand{\dx}{\,{\rm d}x}
\newcommand{\dy}{\,{\rm d}y}

\newcommand{\dt}{\,{\rm d}t}
\newcommand{\rd}{{\rm d}}

\def\LL{\mathrm{L}} %per gli spazi L^p
\def\supp{\mathrm{supp}} %per il supporto

\newcommand{\A}{\mathcal{L}}
\newcommand{\AI}{\mathcal{L}^{-1}}

\newcommand{\n}{F}
\newcommand{\nl}{{F^*}}

\newcommand{\p}{{\delta^\gamma}} %notation for the distance to the gamma, replacing \Phi_1

\newcommand{\ka}{\overline{\kappa}}

\newcommand{\kb}{\underline{\kappa}}

\renewcommand{\k}{\kappa}
\newcommand{\K}{{\mathbb G}}

\newcommand{\RR}{\mathbb{R}}
\newcommand{\NN}{\mathbb{N}}

\newcommand{\B}{\mathcal{B}}

\newcommand{\ve}{\varepsilon}

\def\ee{\mathrm{e}} %per l'esponenziale
\def\dist{\mathrm{dist}} %per la distanza
\def\diam{\mathrm{diam}} %per la distanza
\def\dom{\mathrm{dom}} %per il dominio

\def\qed{\,\unskip\kern 6pt \penalty 500
\raise -2pt\hbox{\vrule \vbox to8pt{\hrule width 6pt
\vfill\hrule}\vrule}\par}

\def\quotient#1#2{\raise1ex\hbox{$#1$}\Big/\lower1ex\hbox{$#2$}}

\definecolor{darkblue}{rgb}{0.05, .05, .65}
\definecolor{darkgreen}{rgb}{0.1, .65, .1}
\definecolor{darkred}{rgb}{0.8,0,0}

\begin{document}
\title{\textbf{Sharp boundary behaviour of solutions\\
to semilinear nonlocal elliptic equations}}

\author{\Large Matteo Bonforte, %$^{\,a,\,b}$ %\footnote{e-mail address:~matteo.bonforte@uam.es}
Alessio Figalli %$^{\,c}$
and Juan Luis V\'azquez %$^{\,a,\,d}$\\} %\footnote{e-mail address:~juanluis.vazquez@uam.es}
}
\date{} %%  this cancels date in article format

\maketitle
%\vspace{-1cm}
\begin{abstract}
We investigate  quantitative properties of nonnegative solutions $u(x)\ge 0$ to the semilinear diffusion equation \ $\A u= f(u)$, posed in a bounded domain $\Omega\subset \RR^N$ with appropriate homogeneous Dirichlet or outer  boundary conditions.  The operator $\A$ may belong to a quite general class of linear operators that include the standard Laplacian, the two most common definitions of the fractional Laplacian $(-\Delta)^s$ ($0<s<1$) in a bounded domain with zero Dirichlet  conditions, and a number of other nonlocal versions. The nonlinearity $f$ is increasing and looks like a power function $f(u)\sim u^p$, with $p\le 1$.

The aim of this paper is to show sharp quantitative boundary estimates based on a new iteration process. We also prove that, in the interior, solutions are H\"older continuous and even classical (when the operator allows for it). In addition, we get H\"older continuity up to the boundary.

Particularly interesting is the behaviour of solution when the number $\frac{2s}{1-p}$ goes below the exponent $\gamma \in(0,1]$ corresponding to
the H\"older regularity of the first eigenfunction $\A\Phi_1=\lambda_1 \Phi_1$.
Indeed a change of boundary regularity happens in the different regimes $\frac{2s}{1-p} \gtreqqless \gamma$,
and in particular a logarithmic correction appears in the ``critical'' case $\frac{2s}{1-p} = \gamma$.

For instance, in the case of the spectral fractional Laplacian, this surprising boundary behaviour appears in the range $0<s\leq (1-p)/2$.

\end{abstract}

\


%\medskip


\vfill
\small
\noindent {\sc Addresses:}

\noindent Matteo Bonforte. Departamento de Matem\'{a}ticas, Universidad
Aut\'{o}noma de Madrid,\\ Campus de Cantoblanco, 28049 Madrid, Spain.
e-mail address:~\texttt{matteo.bonforte@uam.es}
%Web-page:~\texttt{http://www.uam.es/matteo.bonforte

\noindent Alessio Figalli. ETH Z\"urich, Department of Mathematics,
R\"amistrasse 101,\\ 8092 Z\"urich, Switzerland.
E-mail:\texttt{~alessio.figalli@math.ethz.ch}
%Web-page:\texttt{~http://www.ma.utexas.edu/users/figalli/

\noindent Juan Luis V\'azquez. Departamento de Matem\'{a}ticas, Universidad
Aut\'{o}noma de Madrid,\\ Campus de Cantoblanco, 28049 Madrid, Spain.  e-mail address:~\texttt{juanluis.vazquez@uam.es}

\

\

\normalsize

%%%%%%%%%%%%%%%%%%%%%%%%%%%%%%%%%%%%%%%%%%%%%%%%%%%%%%%%%%%%%%%%%%
\section{Introduction}

In this paper we address the question of obtaining a priori estimates, positivity, upper and lower boundary behaviour, Harnack inequalities, and regularity for nonnegative solutions to Semilinear Elliptic Equations of the form\vspace{-1mm}
%%%%%%%%%%%%%%%%%%%%%%%%%%%%%%%%%%%%%
\begin{equation}\label{FPME.equation}
\A\,u=f(u)  \qquad\mbox{posed in } \Omega\,,\vspace{-1mm}
\end{equation}
%%%%%%%%%%%%%%%%%%%%%%%%%%%%%%%%%%%%%
where $\Omega\subset \RR^N$ is a bounded domain with smooth boundary, $N\ge 2$,  $f:\RR \to \RR$ is a monotone nondecreasing function with $f(0)=0$, and $\A$ is a linear operator, possibly of nonlocal type (the basic examples being the fractional Laplacian operators, but the classical Laplacian operator is also included).  Since the problem is posed in a bounded domain we need boundary conditions, or exterior  conditions in the nonlocal case,   that we assume of Dirichlet type and will be included in the functional definition of operator $\A$. This theory covers a quite large class of local and nonlocal operators and nonlinearities.  The operators  $\A$ include the three most common choices of fractional Laplacian operator with Dirichlet conditions but also many other operators  that  are described in Section \ref{sec.hyp.FL}, see  also \cite{BV-PPR2-1, BFV-Parabolic}. In fact, the interest of the theory we develop lies in the wide applicability. The problem is posed in the context of weak dual solutions, which has been proven to be very convenient for the parabolic theory, and is also convenient in the elliptic case.

The focus of the paper is obtaining  a priori estimates and regularity. The a priori estimates are upper bounds for solutions of both signs and lower bounds for nonnegative solutions.    A basic principle in
the paper is that  sharp boundary estimates may   depend not only on $\A$ but also on the behaviour of the nonlinearity $f(u)$ near $u=0$. For this reason we  assume that the nonlinearity $f$ looks like a power with linear or sublinear growth, namely $f(u)\sim u^p$ for some $0<p<1$ when $0\le u\le 1$,   and in that case we identify the range of parameters where the more complicated behaviour happens.

We point out that, for nonnegative solutions, our quantitative inequalities produce sharp behaviour in the interior and near the boundary, both in the case $f(u)=\lambda u$ (the eigenvalue problem) and when $f(u)=u^p$, $p<1$ (the sublinear problem). Our upper and lower bounds will be formulated in terms of the first eigenfunction $\Phi_1$ of $\A$, that under our assumptions will behave like $\Phi_1\asymp  \dist(\cdot, \partial\Omega)^\gamma$ for a certain characteristic power $\gamma\in (0,1]$\,, cf. Section \ref{sec.hyp.FL}. This constant $\gamma$ plays a big role in the theory.

Apart from its own interest, the motivation for this paper comes from companion papers, \cite{BV-PPR2-1, BFV-Parabolic}. In \cite{BV-PPR2-1} a theory for a general class of nonnegative very weak solutions of the parabolic equation $\partial_t u+\A\n(u)=0$ is built, while in \cite{BFV-Parabolic} we address the parabolic regularity theory: positivity, sharp boundary behaviour, Harnack inequalities, sharp H\"older continuity and higher regularity. The proof of such parabolic results relies in part on the elliptic counterparts contained in this paper.\vspace{-1mm}

In this paper we concentrate the efforts in the study of the sublinear case $p\le 1$, since we are motivated by the study of the Porous Medium Equation of the companion paper \cite{BFV-Parabolic}, see also Subsection \ref{ssec.m-p-comp}. The boundary behaviour when $p> 1$  is indeed the same as for $p=1$.\vspace{-1mm}

\noindent\textit{Notation. }Let us indicate here some notation of general use. The symbol $\infty$ will always denote $+\infty$. We also use the notation $a\asymp b$ whenever there exist constants $c_0,c_1>0$ such that $c_0\,b\le a\le c_1 b$\,. We use the symbols $a\vee b=\max\{a,b\}$ and $a\wedge b=\min\{a,b\}$.
We will always consider bounded domains $\Omega$ with smooth boundary, at least $C^{1,1}$. The question of possible lower regularity of the boundary is not addressed here.

%%%%%%%%%%%%%%%%%%%%%%%%%%%%%%%%%%%%%%%%%%%%%%%%%%%%%%%%%%%%%%%%

\section{Basic assumptions and notation}\label{sec.hyp.FL}

In view of the close relation of this study with the parabolic problem,
most of the assumptions on the class of operators $\A$ are the same as in
 \cite{BV-PPR2-1} and \cite{BFV-Parabolic}. We list them for definiteness
and we refer to the references for comments and explanations.

\noindent $\bullet$ {\bf Basic assumptions on $\A$.} The linear operator $\A: \dom(A)\subseteq\LL^1(\Omega)\to\LL^1(\Omega)$ is assumed to be densely defined  and sub-Markovian, more precisely satisfying (A1) and (A2) below:
\begin{enumerate}%[leftmargin=*]\itemsep2pt \parskip3pt \parsep0pt
\item[(A1)] $\A$ is $m$-accretive on $\LL^1(\Omega)$,
\item[(A2)] If $0\le f\le 1$ then $0\le \ee^{-t\A}f\le 1$\,.
\end{enumerate}
The latter can be equivalently written as
\begin{enumerate}
\item[(A2')] If $\beta$ is a maximal monotone graph in $\RR\times\RR$ with $0\in \beta(0)$, $u\in \dom(\A)$\,, $\A u\in \LL^q(\Omega)$\,, $1\le q\le\infty$\,, $v\in \LL^{q/(q-1)}(\Omega)$\,, $v(x)\in \beta(u(x))$ a.e., then
    \[
    \int_\Omega v(x)\,\A u(x)\dx\ge 0\,.
    \]
\end{enumerate}
Such assumptions are the starting   hypotheses   proposed in the paper \cite{BV-PPR2-1} in order to deal with the parabolic problem $\partial_t u+\A F(u)=0$.
Further theory depends on finer properties of the representation kernel of  $\A$, as follows.

\noindent $\bullet$ {\bf Assumptions on $\AI$.} In other to prove our quantitative estimates, we need to be more specific about operator $\A$. Besides satisfying (A1) and (A2), we will assume that it has a left-inverse $\AI: \LL^1(\Omega)\to \LL^1(\Omega)$  with a kernel $\K$ such that
\[
\AI[f](x)=\int_\Omega \K(x,y)f(y)\dy\,,
\]
and that moreover satisfies at least one of the following estimates, for some $s\in (0,1]$:

\noindent - There exists a constant $c_{1,\Omega}>0$ such that for a.e. $x,y\in \Omega$\,:
\[\tag{K1}
0\le \K(x,y)\le c_{1,\Omega}\,|x-y|^{-(N-2s)}\,.
\]
- There exist constants $\gamma\in (0,1]$\,, $c_{0,\Omega},c_{1,\Omega}>0$ such that for a.e. $x,y\in \Omega$\,:
\[\tag{K2}
c_{0,\Omega}\,\p(x)\,\p(y) \le \K(x,y)\le \frac{c_{1,\Omega}}{|x-y|^{N-2s}}
\left(\frac{\p(x)}{|x-y|^\gamma}\wedge 1\right)
\left(\frac{\p(y)}{|x-y|^\gamma}\wedge 1\right)
\]
where we adopt the notation $\p(x):=\dist(x, \partial\Omega)^\gamma$.
Hypothesis (K2) introduces an exponent $\gamma$, which is a  characteristic of the operator and will play a big role in the results.
Notice that defining an inverse operator $\AI$ implies that we are taking into account the Dirichlet boundary conditions.

\medskip

\noindent - The lower bound of assumption (K2) is weaker than the best known estimate on the Green function for many examples under consideration; a stronger inequality holds in many cases:
\[\tag{K4}
\K(x,y)\asymp \frac{1}{|x-y|^{N-2s}}
\left(\frac{\p(x)}{|x-y|^\gamma}\wedge 1\right)
\left(\frac{\p(y)}{|x-y|^\gamma}\wedge 1\right)\,.
\]

\noindent\textbf{The role of the first eigenfunction of $\A$. }Under the assumption (K1) it is possible to show that the operator $\AI$ has a first nonnegative and bounded eigenfunction $0\le \Phi_1\in\LL^\infty(\Omega)$\,, satisfying  $\A\Phi_1=\lambda_1\Phi_1$ for some $\lambda_1>0$, cf. Proposition  \ref{prop.AI.Phi1}. As a consequence of (K2), we show in Proposition \ref {cor.AI.Phi1} that the first eigenfunction satisfies
\begin{equation}\label{Phi1.est}
\Phi_1(x)\asymp \p(x)=\dist(x,\partial\Omega)^\gamma\qquad\mbox{for all }x\in \overline{\Omega}\,,
\end{equation}
hence it encodes the parameter $\gamma$\,, which takes care of describing the boundary behaviour, as first noticed in \cite{BV-PPR1}.

We will also show that all possible eigenfunctions  of $\AI$ satisfy the bound $|\Phi_n|\le \k_n \p$, cf. Proposition \ref{cor.eigenfunctions}.
Recall that we are assuming that the boundary of the domain $\Omega$ is smooth enough, for instance $C^{1,1}$.

In view of \eqref{Phi1.est}, we
can rewrite (K2) and (K4) in the following equivalent forms:  There exist constants $\gamma\in (0,1]$\,, $c_{0,\Omega},c_{1,\Omega}>0$ such that for a.e. $x,y\in \Omega$\,:
\[\tag{K3}
c_{0,\Omega}\Phi_1(x)\Phi_1(y)\le \K(x,y)\le \frac{c_{1,\Omega}}{|x-y|^{N-2s}}
\left(\frac{\Phi_1(x)}{|x-y|^\gamma}\wedge 1\right)
\left(\frac{\Phi_1(y)}{|x-y|^\gamma}\wedge 1\right)\,,
\]
and
\[\tag{K5}
\K(x,y)\asymp \frac{1}{|x-y|^{N-2s}}
\left(\frac{\Phi_1(x)}{|x-y|^\gamma}\wedge 1\right)
\left(\frac{\Phi_1(y)}{|x-y|^\gamma}\wedge 1\right)\,.
\]
We keep the labels (K2), (K4), (K3) and (K5) to be consistent with the papers \cite{BV-PPR2-1,BFV-Parabolic}.\vspace{-1mm}
\subsection{Main Examples}\vspace{-1mm}The theory applies to a number of operators, mainly nonlocal but also  local. We will just list the main cases with some comments, since we have already presented a detailed exposition in \cite{BV-PPR2-1, BFV-Parabolic} that applies here.
In all the examples below, the operators satisfy assumptions $(A1)$ and $(A2)$ and $(K4)$.

As far as fractional Laplacians are concerned, there are at least three different and non-equivalent operators when working on bounded domains, that we  call Restricted Fractional Laplacian (RFL) , the Spectral Fractional Laplacian (SFL) and the Censored Fractional Laplacian (CFL), see Section 3 of \cite{BV-PPR2-1} and Section 2.1 of \cite{BFV-Parabolic}. A good functional setup both for the SFL and the RFL in the framework of fractional Sobolev spaces can be found in \cite{BSV2013}.

For the application  of our results to these cases, it is important to recall that for the RFL $\gamma=s\in (0,1)$, for the CFL $\gamma=s-1/2$ and $s\in (1/2,1)$, while for SFL $\gamma=1$ and $s\in (0,1)$.

There are a number of other operators to which our theory applies: (i) Fractional operators with more general kernels of RFL and CFL type, under some $C^\alpha$ assumptions on the kernel; (ii) Spectral powers of uniformly elliptic operators with $C^1$ coefficients; (iii) Sums of two fractional operators; (iv) Sum of the Laplacian   and a nonlocal operator of L\'evy-type;   (v) Schr\"odinger equations for non-symmetric diffusions; (vi) Gradient perturbation of restricted fractional Laplacians; (vii) Relativistic stable processes, and many other examples more. These examples are presented in detail in Section 3 of \cite{BV-PPR2-1} and Section 10 of \cite{BFV-Parabolic}. Finally, it is worth mentioning that our arguments readily extend to operators on manifolds for which the required bounds on $\K$ hold.

\section{\bf Outline of the paper and main results}\label{sec.MainRes}

In this section we give a overview of the results that we obtain in this paper. Although the first two examples (the linear problem and the eigenvalue problem) are easier and rather standard, some of the results proved in these settings are preparatory for the semilinear problem $\A u=f(u)$, which is the main focus of this paper. In addition, since we could not find a precise reference for (i) and (ii) below in our generality, we present all the details.

\noindent (i) \textit{The linear equation. }We consider the linear problem $\A u=f$ with $0\le f\in \LL^{q'}$ with   $q'=q/(q-1)>N/2s$,   and we show that nonnegative solutions behave at the boundary as follows
\begin{equation}\label{1.4}
\kappa_1\|f\|_{\LL^1_{\Phi_1}(\Omega)} \Phi_1(x)\le  u(x)\le\kappa_2  \|f\|_{\LL^{q'}(\Omega)}\B_q(\Phi_1(x))
\end{equation}
where $\B_q$ is the function defined in \eqref{Lem.Green.est.Upper.II.B}, and depends on the value of $q\in (0, N/(N-2s))$, while $\kappa_1,\kappa_2>0$ depend only on $N,q,\Omega$. See details in Section \ref{sec.lin.f}.

\medskip

\noindent (ii) \textit{Eigenvalue problem. }We  prove a set of a priori estimates for the eigenfunctions, i.e. solutions the Dirichlet problem for the equation $\A \Phi_k=\lambda_k \Phi_k$. We first prove that, under assumption (K1),  eigenfunctions exist and are bounded,
see Proposition \ref{prop.AI.Phi1} and Lemma \ref{Lem.bound.K1}. Then, under assumption (K2), we show the  boundary estimates
\[
\Phi_1(x)\asymp \dist(x,\partial\Omega)^\gamma\qquad\mbox{and}\qquad |\Phi_n(x)|\lesssim \dist(x,\partial\Omega)^\gamma,
\]
see Section \ref{sec.linear.eigen} for more details and results.
Boundary estimates have been proven in various settings,  especially for the common fractional operators (RFL and SFL),   see for instance \cite{BG,BSV2013,Cabre-Tan,CDDS,CS2016,Davies1,Grub1, RosSer, RosSer1, RosSer2, SV1,SV2003}.

\noindent(iii) \textit{Semilinear equations. }
This is the core of the paper, and our main result concerns sharp boundary behaviour.   In Section \ref{sec.semilin} we show that all nonnegative solutions to the semilinear equation \eqref{FPME.equation} with $f(u)\asymp u^p$, $0<p<1$,    satisfy the following sharp estimates whenever $2s+p\gamma \ne \gamma$:
\begin{equation}\label{01}
\kappa_1 \,\dist(x,\partial\Omega)^{\mu} \le
u(x) \le \kappa_2\,  \dist(x,\partial\Omega)^{\mu}\qquad\mbox{for all $x\in \Omega$}\,.
\end{equation}
Here $\kappa_1,\kappa_2>0$ depend only on $N,s,\gamma,p,\Omega$, and the exponent is given by\vspace{-1mm}
\begin{equation}
\label{eq:mu}
\mu:=\gamma\wedge2s/(1-p).\vspace{-1mm}
\end{equation}
Note that $\mu=\gamma$ (i.e., it is independent of $p$) whenever $2s+\gamma p>\gamma$. In particular the ``exceptional value'' $\mu<\gamma$ does not appear neither when $p=1$ or when $s\ge1/2$,   nor in the case of the RFL or CFL.   See also the survey \cite{RosOton1}.  When $2s+\gamma p=\gamma$ (i.e. in the limit case $\mu=\gamma=2s/(1-p)$) a logarithmic correction
appear, and we prove the following sharp estimate:\vspace{-1mm}
\begin{equation}\label{01.log}
u(x) \asymp \dist(x,\partial\Omega)^{\gamma} \left(1+|\log \dist(x,\partial\Omega)|^{\frac{1}{1-p}}\right)\qquad\mbox{for all $x\in \Omega$}\,.\vspace{-1mm}
\end{equation}

\noindent (iv)\textit{ Regularity. }In Section \ref{sec.regularity} we prove that, both in the linear and semilinear case, solutions are H\"older continuous and even classical in the interior (whenever the operator allows it).
In addition we prove that they are H\"older continuous up to the boundary with a sharp exponent. Regularity estimates have been extensively studied: as far as interior H\"older regularity is concerned, see for instance \cite{CDDS, Ka, BCI, CS2, CS3, Si, Tan}; for boundary regularity see \cite{CS2016, Grub1, RosSer, RosSer1, RosSer2, SV1}; for interior  Schauder estimates, see  \cite{BaFiVa, DK}.

\noindent\textbf{Remark. } The results apply without changes in dimension  $N=1$ when $s \in (0,1/2)$.

\noindent {\bf Method and generality. }The usual approach to prove a priori estimates for both linear and semilinear equations, relies De Giorgi-Nash-Moser technique, exploiting energy estimates, and Sobolev and Stroock-Varopoulos inequalities. In addition, extension methods \`a la Caffarelli-Silvestre \cite{Caffarelli-Silvestre} turn out to be very useful. However, due to the generality of the class of operators considered here, such extension is not always possible. Hence, we develop a new approach where we concentrate on the properties of the properties of the Green function of $\A$. In particular, once good linear estimates for the Green function are known, we proceed through a delicate iteration process to establish sharp boundary behaviour of solutions even in a nonlinear setting,   see Propositions \ref{cor.AI.Phi1} and \ref{Prop.Green.2aaa}, and Lemmata \ref{Lem.Green.2}, \ref{Lem.Green.2b}  and \ref{Lem.Green.3}.

\section{The linear problem. Potential and boundary estimates.}\label{sec.lin.f}
In this section we prove estimates on the boundary behaviour of solutions to the linear elliptic problem with zero Dirichlet boundary conditions
\begin{equation}\label{Elliptic.prob.linear.f}
\left\{\begin{array}{lll}
\A u= f &  ~ {\rm in}~  \Omega\\
u=0 & ~\mbox{on the lateral boundary}.
\end{array}
\right.
\end{equation}
The solution to this problem is given by the representation formula
\begin{equation}\label{repr.form.linear}
u(x):=\AI[f](x)=\int_\Omega \K(x,y)f(y)\dy
\end{equation}
whenever $f\in \LL^{q'}$ with $q'>N/2s$. This representation formula is compatible with the concept of weak dual solution that we shall use in the semilinear problem, see Section \ref{sec.semilin}; this can be easily seen by using the definition of weak dual solution and approximating the Green function by means of admissible test functions, analogously to what is done in Subsection \ref{section:proof thms}.
  In the case of SFL and/or of powers of elliptic operators with continuous coefficient, boundary estimates were obtained in \cite{CS2016,CDDS, SV2003}, and for RFL and CFL see \cite{Ka,RosOton1} and references therein. See also Section 3.3 of \cite{BV-PPR2-1} for more examples and references.

The main result of this section is the following theorem.
\begin{thm}\label{Prop.Elliptic.prob.linear.f}
Let $\K$ be the kernel of $\AI$, and assume $(K2)$. Let $u$ be a weak dual solution of the Dirichlet Problem \eqref{Elliptic.prob.linear.f}, corresponding to $0\le f\in \LL^{q'}$ with $q'>N/2s$. Then there exist  positive constants $ \kb_0,\ka_q$ depending on $N,s,\gamma,\Omega,q'$ such that the following estimates hold true
\begin{equation}\label{Elliptic.prob.linear.f.estimates}
\kb_0 \|f\|_{\LL^1_{\Phi_1}(\Omega)} \Phi_1(x)\le u(x)\le\ka_q  \|f\|_{\LL^{q'}(\Omega)}\B_q(\Phi_1(x))\qquad \forall\,x \in \Omega,
\end{equation}
where $q=\frac{q'}{q'-1}\in \left(0,\frac{N}{N-2s}\right)$, and
$\B_q:[0,\infty)\to[0,\infty)$ is defined as follows:
\begin{equation}\label{Lem.Green.est.Upper.II.B}
\B_q(\Phi_1(x_0)):=\left\{\begin{array}{lll}
\Phi_1(x_0)\,, & \qquad\mbox{for }0< q <\frac{N}{N-2s+\gamma}\,,\\[2mm]
\Phi_1(x_0)\,\big(1+ \big|\log\Phi_1(x_0)\big|^{\frac{1}{q}}\big)\,, & \qquad\mbox{for }q = \frac{N}{N-2s+\gamma}\,,\\[2mm]
\Phi_1(x_0)^{\frac{N-q(N-2s)}{q\gamma}}\,, & \qquad\mbox{for }\frac{N}{N-2s+\gamma}<q<\frac{N}{N-2s}\,.\\
\end{array}\right.
\end{equation}
\end{thm}
\noindent\textbf{Remark on the existence of eigenfunctions. }Under assumption (K1) on the kernel $\K$ of $\AI$ we have existence of a positive and bounded eigenfunction $\Phi_1$,  see Subsections  \ref{ssec.existence.eigenfns} and \ref{ssec.boundedeness.eigenfns}; if we further assume (K2) then $\Phi_1\asymp \dist(\cdot\,, \partial\Omega)^\gamma$, cf. Subsection \ref{ssec.bound.beh.eigenfns} for further details.

The proof of the theorem is a simple consequence of the following Lemma
\begin{lem}[Green function estimates I]\label{Lem.Green}Let $\K$ be the kernel of $\AI$, and assume that $(K1)$ holds. Then, for all $0<q<{N}/(N-2s)$, there exist a constant $c_{2,\Omega}(q)>0$ such that
\begin{equation}\label{Lem.Green.est.Upper.I}
\sup_{x_0\in\Omega}\int_{\Omega}\K^q(x , x_0)\dx \le c_{2,\Omega}(q)\,.
\end{equation}
Moreover, if (K2) holds, then  for the same range of $q$ there exists a constant $c_{3,\Omega}(q)>0$ such that,  for all $x_0\in \Omega$,
\begin{equation}\label{Lem.Green.est.Upper.II}
c_{3,\Omega}(q) \Phi_1(x_0)\le \left(\int_{\Omega}\K^q(x , x_0)\dx\right)^{\frac{1}{q}} \le   c_{4,\Omega}(q) \B_q(\Phi_1(x_0))\,,
\end{equation}
where $\B_q:[0,\infty)\to[0,\infty)$ is defined  as in \eqref{Lem.Green.est.Upper.II.B}.
Finally, for all $0\le f\in \LL^1_{\Phi_1}(\Omega)$, $(K2)$ implies that
\begin{equation}\label{Lem.Green.est.Lower.II}
\int_\Omega f(x)\K(x,x_0)\dx\ge c_{0,\Omega} \Phi_1(x_0) \|f\|_{\LL^1_{\Phi_1}(\Omega)}\qquad \text{for all $x_0\in \Omega$}\,.
\end{equation}
\end{lem}
The constants $c_{i,\Omega}(\cdot)$\,, $i=2,3,4,5$\,, depend only on $s,N,\gamma, q, \Omega$, and have an explicit expression given in the proof.

\noindent {\sl Proof of Theorem \ref{Prop.Elliptic.prob.linear.f}.~}
Thanks to \eqref{Lem.Green.est.Upper.I}, the formula
$$
u(x)=\int_\Omega f(y)\K(x,y)\dy
$$
makes sense for $f\in \LL^{q'}(\Omega)$ with $q'>N/2s$.
Now the lower bound is given in \eqref{Lem.Green.est.Lower.II},
while the upper bound follows by \eqref{Lem.Green.est.Upper.II} and H\"older inequality:
\[
u(x)=\int_\Omega f(y)\K(x,y)\dy\le \|f\|_{\LL^{q'}(\Omega)}\|\K(x,\cdot)\|_{\LL^q(\Omega)}\le c_{4,\Omega}(q) \|f\|_{\LL^{q'}(\Omega)} \B_q(\Phi_1(x_0)).\qquad \text{\qed}
\]
\noindent\textsl{Proof of Lemma \ref{Lem.Green}.~}%\label{app.green0}
We split the proof in three steps.

\noindent$\bullet~$\textsc{Step 1. }\textit{Proof of estimate \eqref{Lem.Green.est.Upper.I}. }As consequence of assumption (K1) we obtain
\[\begin{split}
\sup_{x_0\in\Omega}\int_{\Omega}\K^q(x , x_0)\dx
&\le c_{1,\Omega}\sup_{x_0\in\Omega}\int_{\Omega}\frac{1}{|x-x_0|^{q(N-2s)}}\dx\\
&\le c_{1,\Omega}\sup_{x_0\in\Omega}  \int_{B_{{\rm diam}(\Omega)}(x_0)}\frac{1}{|x-x_0|^{q(N-2s)}}\dx\\
&= c_{1,\Omega}\,\frac{N\omega_N}{N-q(N-2s)}\diam(\Omega)^{N-q(N-2s)}=:c_{2,\Omega}(q)\,,
\end{split}
\]
where we used that $\Omega\subset B_{{\rm diam}(\Omega)}(x_0)$ and the notation
$\omega_N=|B_1|$ (recall that by assumption $q(N-2s)<N$).

\smallskip

\noindent$\bullet~$\textsc{Step 2. }\textit{Proof of estimate \eqref{Lem.Green.est.Upper.II}. }We first prove the lower bound of inequality \eqref{Lem.Green.est.Upper.II}. This follows directly from (K2) (see also the equivalent form (K3)):
\[\begin{split}
\int_{\Omega}\K^q(x , x_0)\dx
&\ge c_{0,\Omega} \Phi_1^q(x_0) \int_{\Omega}\Phi_1^q(x) \dx :=c_{3,\Omega}^q(q)\Phi_1^q(x_0)\,.
\end{split}
\]
We next prove the upper bounds of inequality \eqref{Lem.Green.est.Upper.II}. Let us fix $x_0\in \Omega$, and define
\[
R_0:=\Phi_1(x_0)^{1/\gamma}\le \overline{R}:=\|\Phi_1\|_{\LL^\infty(\Omega)}^{1/\gamma}+\diam(\Omega)\,,
\]
so that for any $x_0\in\Omega$ we have $\Omega\subseteq B_{\overline{R}}(x_0)$\,.

Notice that it is not restrictive to assume $0\le \Phi_1(x_0)\le 1$, since we are focusing here on the boundary behaviour, i.e. when $\dist(x, \partial\Omega)\ll 1$ (note that when $\Phi_1(x_0)\ge 1$ we already have estimates \eqref{Lem.Green.est.Upper.I}).\\
Recall now the upper part of (K2) estimates, that can be rewritten in the form
\begin{equation}\label{typeII.Green.est.b}
\K(x,x_0)\le
\frac{c_{1,\Omega}}{|x-x_0|^{d-2s}}%\left(\frac{\Phi_1(x)}{|x-x_0|^\gamma}\wedge 1\right)
\left\{\begin{array}{cl}
\dfrac{\Phi_1(x_0)}{|x-x_0|^\gamma} &\quad\mbox{for any }x\in B_{\overline{R}}(x_0)\setminus B_{R_0}(x_0)\\
1&\quad\mbox{for any }x\in B_{R_0}(x_0)\\
\end{array}\right.
\end{equation}
so that
\[\begin{split}
\int_{\Omega}\K^q(x , x_0)\dx
&\le c_{1,\Omega}^q\left(\int_{B_{R_0}(x_0)}\frac{1}{|x-x_0|^{q(N-2s)}}\dx +\int_{B_{\overline{R}}(x_0)\setminus B_{R_0}(x_0)}\frac{\Phi_1(x_0)^q}{|x-x_0|^{q(N-2s+\gamma)}}\dx\right)\\
&\le c_{1,\Omega}^q\,\omega_N\left(\int_0^{R_0}\frac{r^{N-1}}{r^{q(N-2s)}}\rd r +\Phi_1(x_0)^q\int_{R_0}^{\overline{R}}\frac{r^{N-1}}{r^{q(N-2s+\gamma)}}\rd r\right)
:=(A).
\end{split}
\]
We consider  three cases, depending whether $N-q(N-2s+\gamma)$ is positive, negative, or zero.

- We first analyze the case when $N-q(N-2s+\gamma)>0$\,; recalling that $R_0:=\Phi_1(x_0)^{1/\gamma}$,we have
\[\begin{split}
(A)\le c_{1,\Omega}^q\,\omega_N\left[\frac{R_0^{N-q(N-2s)-q\gamma}}{N-q(N-2s)}
+ \frac{\overline{R}^{N-q(N-2s+\gamma)}}{N-q(N-2s+\gamma)}\right]R_0^{q\gamma}
\leq c_{4,\Omega}^q \Phi_1(x_0)^q.
\end{split}
\]
- Next we analyze the case when $N-q(N-2s+\gamma)<0$\,; using again that $R_0:=\Phi_1(x_0)^{1/\gamma}$, we get
\[\begin{split}
(A)&\le c_{1,\Omega}^q\,\omega_N\left[\frac{1}{N-q(N-2s)}
+ \frac{1}{q(N-2s+\gamma)-N}\right]R_0^{N-q(N-2s)}
= c_{4,\Omega}^q \Phi_1(x_0)^{\frac{N-q(N-2s)}{\gamma}}.
\end{split}
\]
- Finally we analyze the case when $N-q(N-2s+\gamma)=0$\,; again since $R_0:=\Phi_1(x_0)^{1/\gamma}$, it holds
\[\begin{split}
(A)\leq  c_{4,\Omega}^q |\log R_0|\Phi_1(x_0)^q\\
\end{split}
\]
(note that, since we are assuming $R_0\ll 1$,  $-\log R_0 =|\log R_0|\ge 1$). The proof of the upper bound \eqref{Lem.Green.est.Upper.II} is now complete.

\noindent$\bullet~$\textsc{Step 3. }\textit{Proof of estimates \eqref{Lem.Green.est.Lower.II}}. For all $f\in \LL^1_{\Phi_1}(\Omega)$, and $x_0\in\Omega$, the lower bound in (K2) implies
\[
\int_\Omega f(x)\K(x,x_0)\dx\ge c_{0,\Omega} \Phi_1(x_0)\int_\Omega f(x)\Phi_1(x)\dx = c_{0,\Omega} \Phi_1(x_0) \|f\|_{\LL^1_{\Phi_1}(\Omega)}. \qquad \text{\qed}
\]
%
%%%%%%%%%%%%%%%%%%%%%%%%%%%%%%%%%%%%%%%%%%%%%%%%%%%%%%%%%%%%%%%%%%%%%%%%%%%%%%
\section{The eigenvalue Problem}\label{sec.linear.eigen}
%%%%%%%%%%%%%%%%%%%%%%%%%%%%%%%%%%%%%%%%%%%%%%%%%%%%%%%%%%%%%%%%%%%%%%%%%%%%%%%

In this section we will focus the attention on the eigenvalue problem
\begin{equation}\label{Elliptic.prob.eigenvalues}
\left\{\begin{array}{lll}
\A\Phi_k= \lambda_k \Phi_k &  ~ {\rm in}~  \Omega\\
\Phi_k=0 & ~\mbox{on the lateral boundary}
\end{array}
\right.
\end{equation}
 for $k\in \NN$. It is clear, by standard Spectral  theory, that the eigenelements of $\A$ and $\AI$ are the same. We hence focus our study on the ``dual'' problem for $\AI$.

We are going to prove first that assumption (K1) is sufficient to ensure
that the self-adjoint operator $\AI:\LL^2(\Omega)\to \LL^2(\Omega)$ is compact, hence it possesses a discrete spectrum. Then we show that eigenfunctions are bounded. Finally, as a consequence of the stronger assumption (K2), we will obtain the sharp boundary behaviour of the first positive eigenfunction $\Phi_1\asymp \p =\dist(\cdot\,,\partial\Omega)^\gamma$\,, and also optimal boundary estimates for all the other eigenfunctions, namely we prove also that $|\Phi_n|\le \ka_n \p$\,.

\subsection{Compactness and existence of eigenfunctions.}\label{ssec.existence.eigenfns}
Let $\K$ be the kernel of $\AI$, and assume that $(K1)$ holds. Under this assumption we show that $\AI$ is compact, hence it has a discrete spectrum.
\begin{prop}\label{prop.AI.Phi1}
Assume that $\A$ satisfies  $(A1)$ and $(A2)$\,, and that its inverse $\AI$ satisfies $(K1)$.
Then the operator $\AI:\LL^2(\Omega)\to \LL^2(\Omega)$ is compact. As a consequence, $\AI$ possesses a discrete spectrum, denoted by $(\mu_n, \Phi_n)$, with $\mu_n\to 0^+$ as $n\to \infty$. Moreover, there exists a first eigenfunction $\Phi_1\ge 0$ and a first positive eigenvalue $\mu_1=\lambda_1^{-1}>0$\,, such that
\begin{equation}\label{mu1}
0<\lambda_1=\inf_{u\in \LL^2(\Omega)}\frac{\int_\Omega u^2 \dx}{\int_\Omega u\AI u\dx}=\frac{\int_\Omega\Phi_1^2\dx}{\int_\Omega \Phi_1\AI \Phi_1\dx}.
\end{equation}
As a consequence, the following Poincar\'e inequality holds:
\begin{equation}\label{Poincare.AI}
\lambda_1\int_\Omega u\AI u\dx \le \int_\Omega u^2 \dx \qquad\mbox{for all }u\in\LL^2(\Omega)\,.
\end{equation}
\end{prop}

\noindent {\sl Proof of Proposition \ref{prop.AI.Phi1}.~}The proof is divided in several steps. We first prove that the self-adjoint operator $\AI$ is bounded.

\noindent$\bullet~$\textsc{Step 1. }\textit{Boundedness of $\AI$. }We shall prove the following inequality: there exists a constant $C>0$ such that for all $u\in \LL^2(\Omega)$ we have
\begin{equation}\label{prop.AI.Phi1.Step.1.1}
\int_\Omega |\AI u|^2\dx \le C\int_\Omega |u|^2 \dx\,.
\end{equation}
For this, we have
\begin{equation}\label{prop.AI.Phi1.Step.1.2}\begin{split}
\left\|\AI u\right\|_{\LL^2(\Omega)}^2&=\int_\Omega |\AI u|^2\dx
=\int_\Omega \left|\int_\Omega \K(x,y)u(y)\dy\right|^2\dx\\
&\le \int_\Omega \left(\int_\Omega \K(x,y)|u(y)|\dy\right)^2\dx
\le c_{1,\Omega}^2\int_{\RR^N} \left(\int_{\RR^N} \frac{|u(y)|}{|x-y|^{N-2s}}\dy\right)^2\dx\\
&= c_{1,\Omega}^2\int_{\RR^N} \left[(-\Delta_{\RR^N})^{-s}|u|\right]^2\dx \le c_{1,\Omega}^2\left\|(-\Delta_{|\Omega})^{-s}|u|\right\|_{\LL^2(\Omega)}^2\,.
\end{split}
\end{equation}
The last inequality holds because we know that
\[
\|E(u)\|_{H^{-2s}(\RR^N)}=\left(\int_{\RR^N} \left[(-\Delta_{\RR^N})^{-s}|E(u)|\right]^2\dx\right)^{\frac{1}{2}} = \|u\|_{H^{-2s}(\Omega)}=\left\|(-\Delta_{|\Omega})^{-s}|u|\right\|_{\LL^2(\Omega)},
\]
where $E=(r_\Omega^*)^{-1}$ is the inverse of $r_\Omega^*$,  the transpose isomorphism of the restriction operator
$r_\Omega: H^s(\RR^d)\to H^s(\Omega)$; we recall that  ${\rm Ker}(r_\Omega)=\left\{f\in H^s(\RR^d)\;\big|\; r_\Omega(f)=f_{|\Omega}=0\right\}$\,, so that $r_\Omega$ gives the isomorphism
\[
H^s(\Omega)= \quotient{H^s(\RR^d)}{{\rm Ker}(r_\Omega)}
\]
and
\[
H^{-s}(\Omega)=\left(H^s(\Omega)\right)^*
= \left(\quotient{H^s(\RR^d)}{{\rm Ker}(r_\Omega)}\right)^*=\left\{f\in  H^{-s}(\RR^d) \;\big|\; \supp{f}\subseteq\overline{\Omega} \right\}.
\]
We refer to \cite{LM} for further details, see also Section 7.7 of \cite{BSV2013}.\\
Next we recall the Poincar\'e inequality that holds for the Restricted Fractional Laplacian, cf. \cite{BG, ChSo, SV1} and also \cite{BSV2013,BV-PPR2-1}. For all $f\in \LL^2(\Omega)$ such that $(-\Delta_{|\Omega})^{s} f\in \LL^2(\Omega)$ we have that there exists a constant $\lambda_{1,s}>0$ such that
\begin{equation}\label{prop.AI.Phi1.Step.1.4}
\lambda_{1,s}^2 \left\| f\right\|_{\LL^2(\Omega)}^2\le \left\|(-\Delta_{|\Omega})^{s} f\right\|_{\LL^2(\Omega)}^2 \,.
\end{equation}
We apply the above inequality to $f=(-\Delta_{|\Omega})^{-s}|u|$\,, to get
\begin{equation}\label{prop.AI.Phi1.Step.1.5}
\lambda_{1,s}^2 \left\| (-\Delta_{|\Omega})^{-s}|u|\right\|_{\LL^2(\Omega)}^2\le \left\| u\right\|_{\LL^2(\Omega)}^2.
\end{equation}
Combining inequalities \eqref{prop.AI.Phi1.Step.1.2} and \eqref{prop.AI.Phi1.Step.1.5} we obtain \eqref{prop.AI.Phi1.Step.1.1} with $C=\lambda_{1,s}^{-2}c_{1,\Omega}^2>0$\,.

\noindent$\bullet~$\textsc{Step 2. }\textit{The Rayleigh quotient is bounded below: Poincar\'e inequality. }We can compute
\begin{equation}\label{prop.AI.Phi1.Step.1.5b}
\int_\Omega u\AI u\dx\le \left\|\AI u\right\|_{\LL^2(\Omega)} \left\|  u\right\|_{\LL^2(\Omega)} \le C  \left\|  u\right\|_{\LL^2(\Omega)}^2\qquad\mbox{for all $u\in\LL^2(\Omega)$\,,}
\end{equation}
where we have used Cauchy-Schwartz inequality and inequality \eqref{prop.AI.Phi1.Step.1.1} of Step 1.
The above inequality clearly implies that $\lambda_1\ge 1/C>0$, and also proves the Poincar\'e inequality \eqref{Poincare.AI}\,.

\noindent$\bullet~$\textsc{Step 3. }\textit{Compactness. }
Fix $\ve>0$ small and set
$$
\K(x,y)=\K_\ve^1(x,y)+\K_\ve^2(x,y),
$$
where
$$
\K_\ve^1(x,y):=\K(x,y)\chi_{|x-y|\leq \ve},\qquad
\K_\ve^2(x,y):=\K(x,y)\chi_{|x-y|> \ve}.
$$
Note that, by (K1), we can bound
$$
0\leq \K_\ve^1(x,y)\leq G_\ve(x-y),\qquad
G_\ve(z):=c_{1,\Omega}|z|^{-(N-2s)}\chi_{|z|\leq \ve}.
$$
Thus,
for all $u\in \LL^2(\Omega)$ and all $h \in \mathbb R^N$ we have
\begin{equation}\begin{split}\label{step4.Phi.n}
\|\AI u(\cdot+h)-\AI u \|_{\LL^2(\Omega)}&
=\biggl\|\int_{\Omega}\big(\K(\cdot+h,y)-\K(\cdot,y)\big) u(y)\dy\biggr\|_{\LL^2(\Omega)}\\
&\leq \biggl\|\int_{\Omega}\big(\K_\ve^2(\cdot+h,y)-\K_\ve^2(\cdot,y)\big) u(y)\dy\biggr\|_{\LL^2(\Omega)}\\
&\qquad
+\|G_\ve\ast |u|(\cdot +h)\|_{\LL^2(\Omega)}
+\|G_\ve\ast |u|(\cdot)\|_{\LL^2(\Omega)}.
\end{split}
\end{equation}
Now, by Young's convolution inequality, the last two terms can be bounded by
$$
\|u\|_{\LL^2(\Omega)}\|G_\ve(\cdot +h)\|_{\LL^1(\mathbb R^N))}
+\|u\|_{\LL^2(\Omega)}\|G_\ve(\cdot)\|_{\LL^1(\mathbb R^N))}\leq C\|u\|_{\LL^2(\Omega)}\ve^{2s}.
$$
Also, by H\"older inequality we have
\begin{align*}
\biggl\|\int_{\Omega}\big(\K_\ve^2(\cdot+h,y)-\K_\ve^2(\cdot,y)\big) u(y)\dy\biggr\|_{\LL^2(\Omega)}&
\leq \int_{\Omega}\big\|\K_\ve^2(\cdot+h,y)-\K_\ve^2(\cdot,y)\big\|_{\LL^2(\Omega)} |u(y)|\dy\\
&\leq \big\|\K_\ve^2(\cdot+h,\cdot)-\K_\ve^2(\cdot,\cdot)\big\|_{\LL^2(\Omega\times \Omega)} \|u\|_{\LL^2(\Omega)}
\end{align*}
Note that, because of $(K1)$,
 $\K_\ve^2$ is bounded and therefore it belongs to $L^2$. Thus, for $\ve>0$ fixed, it holds that $\big\|\K_\ve^2(\cdot+h,\cdot)-\K_\ve^2(\cdot,\cdot)\big\|_{\LL^2(\Omega\times \Omega)}\to 0$
as $|h|\to 0$, therefore
\[
\lim_{|h|\to 0}\sup_{\|u\|_{\LL^2(\Omega)} \leq 1}\biggl\|\int_{\Omega}\big(\K_\ve^2(\cdot+h,y)-\K_\ve^2(\cdot,y)\big) u(y)\dy\biggr\|_{\LL^2(\Omega)}=0.
\]
Recalling \eqref{step4.Phi.n},
this proves that
\[
\limsup_{|h|\to 0}\sup_{\|u\|_{\LL^2}(\Omega) \leq 1}\|\AI u(\cdot+h)-\AI u \|_{\LL^2(\Omega)}\leq C\ve^{2s},
\]
and since $\ve>0$ is arbitrary we obtain
\[
\lim_{|h|\to 0}\sup_{\|u\|_{\LL^2}(\Omega) \leq 1}\|\AI u(\cdot+h)-\AI u \|_{\LL^2(\Omega)}=0
\]
Since $\AI$ is linear,
 thanks Riesz-Fr\'echet-Kolmogorov Theorem
 we have proved that the image of any ball in $\LL^2(\Omega)$
 is compact in $\LL^2(\Omega)$ with respect to the strong topology. Hence the operator $\AI$ is compact and has a discrete spectrum.

\noindent$\bullet~$\textsc{Step 4. }\textit{The first eigenfunction and the Poincar\'e inequality. } The first eigenfunction exists in view of the previous step. Finally, the minimality property \eqref{mu1} follows by standard arguments and implies both the non-negativity of $\Phi_1$ and the Poincar\'e inequality \eqref{Poincare.AI}.\qed

\subsection{Boundedness of eigenfunctions.}\label{ssec.boundedeness.eigenfns}

We now show that, under the only assumption (K1), all the eigenfunctions  are bounded,
namely there exists $\overline{K}_n>0$ depending only on $N,s,\Omega$ and $n$, such that
\begin{equation}\label{bounds.Phi.n}
\|\Phi_n\|_{\LL^\infty(\Omega)}\le \overline{K}_n\,.
\end{equation}
Recall that we are considering eigenfunctions $\Phi_n$ normalized in $\LL^2(\Omega)$. The key point to obtain such bounds is that the absolute value of eigenfunctions satisfies an integral inequality:
\begin{equation}\label{Phi.k.subsol}
\big|\Phi_n(x_0)\big|=\lambda_n\big|\AI\Phi_n(x_0)\big|=\lambda_n\left|\int_{\Omega} \Phi_n(x)\K(x,x_0)\dx\right|
\le \lambda_n \int_{\Omega} \big|\Phi_n(x)\big|\, \K(x,x_0)\dx.
\end{equation}
Thus $0\le u(x):=\big|\Phi_n(x)\big|$ satisfies the hypothesis of Lemma \ref{Lem.bound.K1} below.

\begin{lem}\label{Lem.bound.K1}Assume that $\A$ satisfies  $(A1)$ and $(A2)$, and that its inverse $\AI$ satisfies $(K1)$. If $u\in \LL^2(\Omega)$ is nonnegative and satisfies
\begin{equation}\label{Lem.Green.lineare.hyp-upper.K1}
 u(x_0) \le \kappa_0\int_{\Omega} u(x)\K(x,x_0)\dx\,,
\end{equation}
then there exists a constant $\ka>0$ such that the following sharp upper bound holds true:
\begin{equation}\label{Lem.Green.lineare.upper.K1}
\|u\|_{L^\infty(\Omega)}
\le \ka \|u\|_{\LL^{1}(\Omega)}\,,
\end{equation}
\end{lem}
\noindent {\sl Proof.~}
The boundedness follows by the Hardy-Littlewood-Sobolev (HLS) inequality, through a finite iteration. The HLS reads
$$
\left\|(-\Delta_{\RR^N})^{-s}f\right\|_{\LL^q(\RR^N)}  \le \mathcal{S}_p \|f\|_{\LL^p(\RR^N)}\,,\qquad\mbox{for all $0<s<N/2$ and }q=\frac{pN}{N-2s p}>p>1,
$$
see \cite{LM} or \cite{BSV2013} and references therein. We will use HLS in the following iterative form:
\begin{equation}\label{prop.AI.Phi1.Step.4.2}
\left\|(-\Delta_{\RR^N})^{-s}f\right\|_{\LL^{p_{k+1}}(\RR^N)}  \le \mathcal{S}_{p_k} \|f\|_{\LL^{p_k}(\RR^N)}\,,
\;\mbox{with } p_{k+1}:=\frac{Np_k}{N-2s p_k}=\frac{Np_0}{N-2s(k+1)p_0}\,.
\end{equation}
Indeed, for all $f\ge 0$\,, as a consequence of (K1) we have that
$$
\AI f(x)=\int_\Omega f(y)\K(x,y)\dy \le c_{1,\Omega}\int_{\RR^N} \frac{f(y)}{|x-y|^{N-2s}}\dy=(-\Delta_{\RR^N})^{-s}f(x)
$$
so that
\begin{equation}\label{prop.AI.Phi1.Step.4.4}
\left\|(-\Delta_{\RR^N})^{-s}f\right\|_{\LL^q(\Omega)}\le c_{1,\Omega}^{1/q}\left\|(-\Delta_{\RR^N})^{-s}f\right\|_{\LL^q(\RR^N)}.
\end{equation}
As a consequence of  \eqref{prop.AI.Phi1.Step.4.2} and \eqref{prop.AI.Phi1.Step.4.4}\,, we obtain the following inequality for all $f\ge 0$ supported in $\overline{\Omega}$:
\begin{equation}\label{prop.AI.Phi1.Step.4.5}
\left\|\AI f\right\|_{\LL^{p_{k+1}}(\Omega)}\le c_{1,\Omega}^{1/p_{k+1}}\left\|(-\Delta_{\RR^N})^{-s}f\right\|_{\LL^{p_{k+1}}(\RR^N)}\le \mathcal{S}_k \|f\|_{\LL^{p_k}(\Omega)}\,,
\end{equation}
with $p_k$ as in \eqref{prop.AI.Phi1.Step.4.2} and $\mathcal{S}_k=c_{1,\Omega}^{1/p_{k+1}}\mathcal{S}_{p_k}$\,. \\
Now, using the inequality \eqref{Lem.Green.lineare.hyp-upper.K1} satisfied by $ u$\,, namely $u\le \mu_1^{-1}\AI u$ and  setting $p_0:=2$ and $\LL^{p_k}:=\LL^{p_k}(\Omega)$, we get
\begin{equation}\label{prop.AI.Phi1.Step.4.6}
\begin{split}
\| u\|_{\LL^{p_{k+1}} }&\le\mu_1\|\AI u\|_{\LL^{p_{k+1}} }
\le \mu_1\mathcal{S}_k\| u\|_{\LL^{p_k} }\le \mu_1^2\mathcal{S}_k\|\AI u\|_{\LL^{p_k} }
\le \mu_1^2\mathcal{S}_k\mathcal S_{k-1}\| u\|_{\LL^{p_{k-1}} }\\
&\le\mu_1^3\mathcal{S}_k\mathcal S_{k-1}\|\AI u\|_{\LL^{p_{k-1}} }
\le\ldots\le \mu_1^{k+1}\biggl(\prod_{j=0}^k\mathcal{S}_j\biggr)\| u\|_{\LL^{2} }.
\end{split}
\end{equation}
After a finite number of steps, namely until
\[
k\ge \frac{N}{2s p_0}-2,\qquad\mbox{we then have}\qquad p_{k+1}>\frac{N}{2s}\,.
\]
Thus, thanks to \eqref{prop.AI.Phi1.Step.4.6},
we get
\[
0\le \lambda_1^{-1} u(x)=\int_\Omega u(y)\K(x,y)\dy \le \| u\|_{\LL^{p_{k+1}}(\Omega)}\|\K(x,\cdot)\|_{\LL^{p'_{k+1}}(\Omega)}\le c_{2,\Omega}(p'_{k+1})\ka\| u\|_{\LL^{2}}\,.
\]
where we have used inequality \eqref{Lem.Green.est.Upper.I} with $1< p'_{k+1}< N/(N-2s)$ (because  $p'_{k+1}=\frac{p_{k+1}}{p_{k+1}-1}$ and $p_{k+1}>N/2s$), so that the constant $c_{2,\Omega}(p'_{k+1})$ is finite.

Hence, this proves that
$$
\|u\|_{L^\infty(\Omega)}
\le \hat\kappa \|u\|_{\LL^{2}(\Omega)}.
$$
Now,
it follows by H\"older inequality that
$$
\|u\|_{L^\infty(\Omega)}
\le \hat\kappa \|u\|_{\LL^{2}(\Omega)}
\leq \hat\kappa \|u\|_{\LL^{\infty}(\Omega)}^{1/2} \|u\|_{\LL^{1}(\Omega)}^{1/2},
$$
thus
$$
\|u\|_{L^\infty(\Omega)}
\le \hat\kappa^2 \|u\|_{\LL^{1}(\Omega)},
$$
concluding the proof.\qed

%%%%%%%%%%%%%%%%%%%%%%%%%%%%%%%%%%%%%%%%%%%%%%%%%%%%%%%%%%%%%%%%%%%%%%%%%%%%%%%%

\subsection{Boundary behaviour of eigenfunctions}\label{ssec.bound.beh.eigenfns}

\begin{prop}\label{cor.AI.Phi1}
Assume that $\A$ satisfies  $(A1)$ and $(A2)$\,, and that its inverse $\AI$ satisfies $(K2$). Then the first eigenfunction satisfies the boundary estimates \eqref{Phi1.est}, namely there exist two positive constants $\ka,\kb$ depending only on $N,s,\Omega$ such that, for all $x\in \overline{\Omega}$,
\[
\kb\, \dist(x,\partial\Omega)^\gamma \le \Phi_1(x) \le \ka\, \dist(x,\partial\Omega)^\gamma\,.
\]
\end{prop}
Because the proof of this result is a modification of the one of Lemma \ref{Lem.Green.2b}
in the case $F(u)=\lambda u$, we postpone it after the proof of Lemma \ref{Lem.Green.2b}, at the end of Subsection \ref{app.green1}.

\begin{prop}\label{cor.eigenfunctions}
Assume that $\A$ satisfies  $(A1)$ and $(A2)$\,, and that its inverse $\AI$ satisfies $(K2)$.
Let $\Phi_n$ be the $n$-th eigenfunction of $\AI$ with unitary $\LL^2$ norm.
Then there exist a positive constant $\ka_n$ depending only on $N,s,n,\Omega$ such that for all $x\in \overline{\Omega}$
\begin{equation}\label{eigenfunctions.bounds}
|\Phi_n(x)| \le \ka_n\, \dist(x,\partial\Omega)^\gamma\,.
\end{equation}
\end{prop}

To prove this second result, we first state and prove
a general Lemma about sub and super solutions

\begin{lem}[Green function estimates III]\label{Lem.Green.lineare}
Assume that $\A$ satisfies  $(A1)$ and $(A2)$\,, and that its inverse $\AI$ satisfies $(K2)$. Let $u:\Omega \to [0,\infty)$ be an integrable function. \\
\noindent\textsc{Upper Bounds. }Assume that for a.e. $x_0\in\Omega$
\begin{equation}\label{Lem.Green.lineare.hyp-upper}
 u(x_0) \le \kappa_0\int_{\Omega} u(x)\K(x,x_0)\dx\,.
\end{equation}
Then, there exists a constant $\ka>0$ such that the following sharp upper bound holds true for a.e. $x_0\in \Omega$\,:
\begin{equation}\label{Lem.Green.lineare.upper}
u(x_0) \le \ka \|u\|_{\LL^1(\Omega)}\dist(x_0,\partial\Omega)^\gamma\,.
\end{equation}
\textsc{Lower bounds. }Assume that for all $x_0\in\Omega$
\begin{equation}\label{Lem.Green.lineare.hyp-lower}
\n(u(x_0))\ge \kappa_0\int_{\Omega} u(x)\K(x,x_0)\dx\,.
\end{equation}
Then, we have the following sharp lower bound for a.e. $x_0\in \Omega$:
\begin{equation}\label{Lem.Green.lineare.lower}
u(x_0) \ge \kb\|u\|_{\LL^1_{\p}(\Omega)}\dist(x_0,\partial\Omega)^\gamma\,.
\end{equation}
Here, $\ka \,,\,\kb >0$ only depend on $s,\gamma,  N,\Omega,\kappa_0$\,.
\end{lem}
\noindent {\sl Proof.~}

\noindent$\bullet~$\textsc{Upper bounds. }From Lemma \ref{Lem.bound.K1} we deduce that $u$ is bounded, more precisely that it satisfies the estimates
\[
\|u\|_{\LL^\infty(\Omega)}\le \ka_0 \|u\|_{\LL^1(\Omega)}.
\]
where $\ka_0$ depends only on $N,s,\gamma,\Omega$\,.

Then we consider the function $\hat u=u/\|u\|_{\LL^1(\Omega)}$\,, which has unitary $\LL^1$ norm.
Applying the argument used to prove the upper bounds of Proposition \ref{cor.AI.Phi1} to the function $\hat u$ in place of $\Phi_1$, one  concludes that
\[
\hat u(x_0) \le \ka \,\dist(x_0,\partial\Omega)^\gamma,\qquad\mbox{that is}\qquad u(x_0) \le \ka \|u\|_{\LL^1(\Omega)} \dist(x_0,\partial\Omega)^\gamma,
\]
which proves \eqref{Lem.Green.lineare.upper}.

\noindent$\bullet~$\textsc{Lower bounds. }The proof is similar to the one of estimates \eqref{Lem.Green.est.Lower.II}, cf. Step 3 in the proof of Lemma \ref{Lem.Green}:
\[
\int_\Omega f(x)\K(x,x_0)\dx\ge c_{0,\Omega} \p(x_0)\int_\Omega f(x)\p(x)\dx
= c_{0,\Omega} \p(x_0) \|f\|_{\LL^1_{\p}(\Omega)}.\qquad \text{\qed}
\]
\medskip

\noindent {\sl Proof of Proposition \ref{cor.eigenfunctions}.~}Since the eigenfunction $\Phi_n$ satisfies the equation $\lambda_n^{-1}\Phi_n=\AI \Phi_n$, we have already shown in \eqref{Phi.k.subsol} that $|\Phi_n|\le \lambda_n \AI\big(|\Phi_n|\big)$,
so that $0\le u(x):=\big|\Phi_n(x)\big|$ satisfies the hypothesis of Lemma \ref{Lem.Green.lineare}\,, hence \eqref{Lem.Green.lineare.upper} hold true. Since $c_\Omega\|\Phi_n\|_{\LL^1(\Omega)}\leq \|\Phi_n\|_{\LL^2(\Omega)}=1$, the result follows.\qed

%%%%%%%%%%%%%%%%%%%%%%%%%%%%%%%%%%%%%%%%%%%%%%%%%%%%%%%%%%%%%%%%%%%%%%%%%%%%%%%%%%%%%%

\section{The Semilinear Problem}\label{sec.semilin}
We consider the following homogeneous Dirichlet problem
\begin{equation}\label{Elliptic.prob}
\left\{\begin{array}{lll}
\A u= f(u) &  ~ {\rm in}~  \Omega\\
u \geq 0 &  ~ {\rm in}~  \Omega\\
u=0 & ~\mbox{on the lateral boundary}
\end{array}
\right.
\end{equation}
where $f:[0,\infty)\to[0,\infty)$\, is concave and increasing with $f(0)=0$. The prototype example is $f(u)=u^p$ with $p\in (0,1)$, as we shall comment later in Subsection \ref{ssec.m-p-comp}.

We begin by defining the concept of weak dual solution for this problem.
Note that, by concavity of $f$ we have $|f(a)|\leq C(1+|a|)$ for all $a\in \mathbb R$, so $f(u)$ is integrable whenever $u$ is so.

\begin{defn}\label{Def.Very.Weak.Sol.Dual.elliptic} A function  $u\ge 0$   is a {\sl weak dual} solution to Problem \eqref{Elliptic.prob} if $u\in \LL^1_{\Phi_1}(\Omega)$ and the identity
\begin{equation}
\int_{\Omega}\AI [f(u)] \,\psi\,\dx
= \int_{\Omega} u\,\psi\,\dx.
\end{equation}
holds for every test function $\psi$ such that $\psi/\Phi_1\in \LL^\infty(\Omega)$\,.
\end{defn}
As far as we know, this definition of weak solution has been introduced by H. Brezis in the case of the classical Laplacian $\A=(-\Delta)$, see also \cite{DR-JFA}. In the parabolic framework, the analogous definition of weak dual solution has been given in \cite{BV-PPR1}, see also \cite{BV-PPR2-1,BSV2013}, and also has been adapted to the elliptic framework there\,.
This concept of weak solution is more general that the usual one (namely, weak solutions used in variational methods, cf.
\cite{BK92,CaSi1, CaSi2, CDDS, Ka, Tan});   indeed, it is not difficult to show that weak solutions are indeed weak dual solutions. Notice that when $f$ does not depend on $u$, this concept of solution is nothing else but the integral representation \eqref{repr.form.linear} of the solution by means of the Green function of $\A$.

The existence of weak solutions follow from standard methods in calculus of variations, with minor modifications. In addition, when $f(u)=u^p$ with $p\in (0,1)$, various techniques to prove uniqueness can be found in the Appendix of \cite{BK92} and they can easily be adapted to  our case. As a consequence, the following result holds:

\begin{thm}[Existence and uniqueness]
Let $f:[0,\infty)\to[0,\infty)$ be a concave increasing function with $f(0)=0$.
Then the semilinear Dirichlet problem \eqref{Elliptic.prob} admits a nontrivial weak dual solution. When $f(u)=u^p$ with $p<1$, then this solution is unique.
\end{thm}
Let us mention that another proof of uniqueness can be done by means of parabolic techniques and separation of variables solutions, see for instance
\cite{AP,JLVmonats} for the classical case $\A=-\Delta$ and \cite{BSV2013, BV-PPR2-1} for more general operators.

\subsection{A priori estimates. Statement of results}

We first show that (sub)solutions are bounded and have a certain decay at the boundary, as a consequence of (K2); this upper boundary behaviour is already sharp in the case $2s>\gamma$, as we discuss later.
\begin{thm}\label{Thm.upper.ell}Assume that $\A$ satisfies  $(A1)$ and $(A2)$\,, and that its inverse $\AI$ satisfies $(K2)$. Let $\n=f^{-1}$ be a positive convex function with $\n(0)=0$. Let $u$ be a weak dual (sub)solution to the Dirichlet Problem \eqref{Elliptic.prob}.
Then, there exists a constant $\kappa_1>0$, depending on $s,N,\gamma,\Omega$ only,  such that for all $x_0\in \Omega$ the following absolute upper bound holds:
\begin{equation}\label{Lem.Green.2.est.Upper.II}u(x_0) \le \kappa_1\B_1(\Phi_1(x_0))\,,
\end{equation}
where $\B_1$ is defined in \eqref{Lem.Green.est.Upper.II.B}.
\end{thm}\vspace{-2mm}
\begin{thm}[Global Harnack Principle]\label{Thm.Elliptic.Harnack.m}Assume that $\A$ satisfies  $(A1)$ and $(A2)$\,, and that its inverse $\AI$ satisfies $(K2)$.
Let $u\ge~ 0$ be a weak dual solution to the Dirichlet Problem \eqref{Elliptic.prob}, where $f$ is a nonnegative increasing function with $f(0)=0$ such that $\n=f^{-1}$ is convex and $\n(a)\asymp a^{1/p}$ when $0\leq a\le 1$, for some $0<p<1$.
Then there exist positive constants $\kb,\ka$, depending only on $\Omega$, $N, p, s,\gamma$,
such that the following sharp absolute bounds hold true: \\
(i) When $2s+p\gamma \ne \gamma$ we have that, for all $x\in \Omega$,
\begin{equation}\label{Thm.Elliptic.Harnack.ineq.m.1}
\kb\, \Phi_1^{\sigma}(x) \le u(x)\le \ka\,\Phi_1^{\sigma}(x) \,,
\end{equation}
where
\begin{equation}
\label{eq:sigma}
\sigma:=1\wedge\frac{2s}{\gamma(1-p)}=\frac{\mu}{\gamma},
\end{equation}
and $\mu$ is as in \eqref{eq:mu}.
When $\sigma<1$\,, the lower bound requires the extra hypothesis $(K4)$\,, otherwise it holds in a non-sharp form with $\sigma=1$\,.\\
(ii) When $2s+\gamma p=\gamma$\,, assuming $(K4)$ we have that, for all $x\in \Omega$,
\begin{equation}\label{Thm.Elliptic.Harnack.ineq.m.log}
\kb\, \Phi_1(x)\left(1+|\log\Phi_1(x) |^{\frac{1}{1-p}}\right) \le u(x)\le \ka\,\Phi_1(x)\left(1+|\log\Phi_1(x) |^{\frac{1}{1-p}}\right)\,.
\end{equation}
If (K4) does not hold, then the upper bound still holds, but the lower bound holds
in a non-sharp form without the extra logarithmic term.
\end{thm}
\noindent\textbf{Remarks. }(i) Here we can appreciate the strong difference between different operators $\A$ reflected in their exponent $\gamma$: recalling that $\Phi_1\asymp\p=\dist(\cdot,\partial\Omega)^\gamma$, we can observe that solutions corresponding to different operators have different boundary behaviour (as we expect for the linear case by the Green function estimates), but we can also appreciate the nontrivial interplay with the nonlinearity, that we shall try to clarify below by means of concrete and relevant examples.

\noindent(iii) As a consequence of Theorem \ref{Thm.Elliptic.Harnack.m}, a more standard form of the \textit{local Harnack inequality }holds: for any $0<R\le \dist(x_0, \partial\Omega)/2$,
\begin{equation}\label{local.Harnack}
\sup_{x\in B_R(x_0)}u(x)\le \k \inf_{x\in B_R(x_0)}u(x)\,.
\end{equation}

We now comment on our result in the case
of the Restricted (RFL), the Spectral (SFL), and the Censored (CFL) Fractional Laplacian, see Section 3 of \cite{BV-PPR2-1} and Section 2.1 of \cite{BFV-Parabolic} for their definition.

\noindent\textbf{The RFL case. }In this case we have always $\gamma=s<2s$, therefore for all $0<s \le 1$ we have
\[
u(x)\asymp\dist(x,\partial\Omega)^s\qquad\qquad\mbox{for   all $x\in \Omega$\,.}
\]
so that $u$ always behaves like the distance to the boundary at the power $s$, the same boundary behaviour of the first eigenfunction $\Phi_1$. The sharp boundary behaviour for eigenfunctions and for the linear equation has been obtained in \cite{RosSer}, see also \cite{Grub1}.

\noindent\textbf{The CFL case. }In this case we have always $\gamma=s-1/2<2s$ for all $1/2<s \le 1$, hence
\[
u(x)\asymp\dist(x,\partial\Omega)^{s-\frac{1}{2}}\qquad\qquad\mbox{for   all $x\in \Omega$\,.}
\]
Thus $u$ always behaves like the distance to the boundary at the power $s-1/2$, the same boundary behaviour of the first eigenfunction $\Phi_1$.

Note that, since $u\in \LL^\infty$ (by Theorem \ref{Thm.upper.ell}), then  $f\sim u^p\in \LL^{\infty}(\Omega)$
and the sharp boundary behaviour for the RFL and the CFL were already known in these particular cases, see for instance \cite{BFR,Ka,RosOton1} and the references therein.

The first two examples may suggest that the boundary behaviour should always be given by $\dist(\cdot,\partial\Omega)^\gamma$. This is actually not the case: for spectral-type Laplacians, the behaviour may change for solutions to $\A u=1$ and to $\A u=\lambda u$, see for instance \cite{CDDS, CS2016}.

\noindent\textbf{The SFL case. }In this case we have $\gamma=1$  hence we can have two different sharp boundary behaviours, in two different range of parameters. First, when $s>  \frac{1-p}{2}$ we obtain
\begin{equation}\label{Thm.Elliptic.Harnack.ineq.m.spectral}
u(x)\asymp \dist(x,\partial\Omega) \qquad\qquad\mbox{for  a.e.  $x\in \Omega$\,.}
\end{equation}
On the other hand, when $s<\frac{1-p}{2}$ we have the following sharp bound
\begin{equation}\label{Lem.Green.2.est.Upper.II.b}
u(x)\asymp\dist(x,\partial\Omega)^{\frac{2s}{1-p}}\qquad\qquad\mbox{for a.e.  $x\in \Omega$\,.}
\end{equation}
Finally, when $s=\frac{1-p}{2}$ we have the following sharp bound
\begin{equation}\label{Lem.Green.2.est.Upper.II.log}
u(x)\asymp\dist(x,\partial\Omega)\left(1+|\log\Phi_1(x) |^{\frac{1}{1-p}}\right)\qquad\qquad\mbox{for a.e.  $x\in \Omega$\,.}
\end{equation}

 \noindent\textbf{Remarks. }(i) In the two estimates above we can appreciate the interplay between the ``scaling power'' $2s/(1-p)$ and the ``eigenfunction power'' $1$. The boundary behaviour of the semilinear equation somehow interpolates between the extremal cases $p=0$, i.e. $\A u= 1$, and $p=1$, i.e. $\A u=\lambda u$.\\
(ii) Our result improves the boundary estimates obtained in \cite{CS2016}, Theorem 1.3.
\subsubsection{Comparison with parabolic estimates}\label{ssec.m-p-comp}
We are interested in comparing the elliptic estimates of this paper with the parabolic estimates of the companion paper \cite{BFV-Parabolic}. There, we consider the nonlinear parabolic equation\vspace{-2mm}
$$
v_t+\A v^m=0, \quad m>1,\,
$$
that admits  separate-variables solutions which take the form $v(t,x)=V(x)t^{-1/(m-1)}$, where
$V(x)$ satisfies the elliptic equation\vspace{-2mm}
$$
\A V^m = \frac{1}{m-1} V.
$$
Hence, setting\ $V^m=u$ we see that  $u$ satisfies the equation $\A u =\frac{1}{m-1} u^p$ with $p=1/m<1$, which is a particular case of equation \eqref{Elliptic.prob} (note that, up to multiplying $u$ by a positive constant, we can remove the moltiplicative term $1/(m-1)$ in the right hand side). See also \cite{BFR,BSV2013, BV-PPR1} and the survey \cite{Vaz2014}. \vspace{-2mm}

In view of this, we introduce in what follows the notation $u=F(V)$ in order to move the nonlinearity in the right hand side: letting $F=f^{-1}$, then $V=f(u)$; in this way, we can restate all the results of this paper in terms of $V$. \vspace{-4mm}

\subsection{Proof of the sharp boundary estimates}\label{subsect.Green.bdry}

In order to prove our upper and lower bound on solutions to \eqref{Elliptic.prob},
we first prove a series of upper and lower bound for sub and subsolutions to a semilinear equation.
Hence the function $u$ in the next statements is just a general bounded function, and our Theorems \ref{Thm.upper.ell} and \ref{Thm.Elliptic.Harnack.m} will follow by applying the following estimates to $f(u)$ in place of $u$.\vspace{-1mm}

We begin  by stating the following important estimates:\vspace{-3mm}
\begin{prop}[Green function estimates II]\label{Prop.Green.2aaa}
Let $\K$ be the kernel of $\AI$ satisfying $(K2)$, and let $F$ be a non-negative convex function with $\n(0)=0$. Let $u:\Omega \to [0,\infty)$ be a bounded measurable function.\\
\noindent\textsc{Upper Bounds.} Assume that for all $x_0\in\Omega$\vspace{-2mm}
\begin{equation}\label{Lem.Green.2.hyp.aaa}
\n(u(x_0))\le \kappa_0\int_{\Omega} u(x)\K(x,x_0)\dx\vspace{-2mm}
\end{equation}
Then there exists a constant $c_{\Omega}>0$ such that  the following absolute upper bounds hold true:  let  $\kappa_1=c_{\Omega}\kappa_0\,\n^{-1}\circ\nl(2c_{\Omega}\kappa_0)$, where $F^*$ denotes the Legendre transform of $F$. Then, with $\B_1$ as in \eqref{Lem.Green.est.Upper.II.B}, we have\vspace{-1mm}
\begin{equation}\label{Lem.Green.2.est.Upper.aaa}
\n(u(x_0))\le \kappa_1\B_1(\Phi_1(x_0))\,,\qquad\mbox{for all $x_0\in \Omega$. }
\end{equation}
Assume further that for $0\leq a\le 1 $ we have $F(a)\ge \kb a^m$ for some $\kb>0$ and $m>1$. Then, when $\gamma<2sm/(m-1)$\,, there exist a constant $\kappa_4>0$ such that,\vspace{-1mm}
\begin{equation}\label{Lem.Green.2.est.Upper.II.b1x}
u^{m}(x_0) \le \kappa_4\kappa_0^{\frac{m}{m-1}}\Phi_1(x_0) \,,\qquad\mbox{for all $x_0\in \Omega$. }
\end{equation}
\noindent On the other hand, when $1\ge \gamma\geq  2s m/(m-1)$, there exist $\kappa_5>0$ such that  $ x_0\in \Omega$
\begin{equation}\label{Lem.Green.3.est.Lower.II.bx}
\begin{split}
 u^{m}(x_0)  &\le
\kappa_{5} \kappa_0^{\frac{m}{m-1}}
\left\{
\begin{array}{ll}
\Phi_1(x_0)\big(1+|\log \Phi_1(x_0)|^{\frac{m}{m-1}}\big)&\text{if $\gamma= 2s m/(m-1)$},\\
\Phi_1(x_0)^{\frac{2sm}{(m-1)\gamma}}&\text{if $\gamma> 2s m/(m-1)$}.\\
\end{array}
\right.
\end{split}
\end{equation}
The upper bounds \eqref{Lem.Green.2.est.Upper.II.b1x} and \eqref{Lem.Green.3.est.Lower.II.bx} are sharp. Here, $\kappa_4,\kappa_5>0$ depend on $s,\gamma, m, N, \kb, \Omega$ only.

\noindent\textsc{Lower bounds. }Assume that $F(a)\le \ka a^m$ for some $m>1$ and for all $0\leq a\le 1$\,, and that
\begin{equation}\label{Lem.Green.3.hyp.aaa}
\n(u(x_0))\ge \kappa_0\int_{\Omega} u(x)\K(x,x_0)\dx\qquad\mbox{for a.e. $x_0\in\Omega$.}
\end{equation}
 Then, there exists a constant $\kappa_{6}>0$ such that, for a.e. $x_0\in \Omega$,
\begin{equation}\label{Lem.Green.3.est.Lower.IIx}
 u^m(x_0)  \ge \kappa_6\kappa_0^{\frac{m}{m-1}} \Phi_1(x_0)\,,
\end{equation}
and the estimates are sharp when $\gamma< 2s m/(m-1)$\,.

\noindent On the other hand, when $1\ge \gamma\geq  2s m/(m-1)$, assumption (K4) implies that there exist a constant $\kappa_7>0$ such that, for a.e.  $x_0\in \Omega$,
\begin{equation*}
\begin{split}
  u^m(x_0)
&\ge
\kappa_{7} \kappa_0^{\frac{m}{m-1}}
\left\{
\begin{array}{ll}
\Phi_1(x_0)\big(1+|\log \Phi_1(x_0)|^{\frac{m}{m-1}}\big)&\text{if $\gamma= 2s m/(m-1)$},\\
\Phi_1(x_0)^{\frac{2sm}{(m-1)\gamma}}&\text{if $\gamma> 2s m/(m-1)$},\\
\end{array}
\right.
\end{split}
\end{equation*}
and the estimates are sharp. Here, $\kappa_6,\kappa_7>0$ depend on $s,\gamma, m, N,\ka, \Omega$ only.
\end{prop}
The proof of Proposition \ref{Prop.Green.2aaa} will be split in three Lemmata, namely Lemmata \ref{Lem.Green.2} and  \ref{Lem.Green.2b} for the upper bounds, and Lemma \ref{Lem.Green.3} for the lower bounds. \vspace{-2mm}
\subsubsection{Green function estimates II. Upper bounds}\label{app.green1}
Let us begin with some preliminary comments. When $f\in \LL^\infty(\Omega)$\,, the estimates of Lemma \ref{Lem.Green} for $q=1$  read as follows:
\begin{equation}\label{rem1.1.Green}
c_{0,\Omega}\Phi_1(x_0) \|f\|_{\LL^1_{\Phi_1}(\Omega)}\le \int_\Omega f(x)\K(x,x_0)\dx\le c_{4,\Omega}\|f\|_{\LL^\infty(\Omega)}\B_1(\Phi_1(x_0)),
\end{equation}
with
\begin{equation}\label{rem1.2.Green}
\B_1(\Phi_1(x_0))=\left\{\begin{array}{lll}
\Phi_1(x_0) & \qquad\mbox{for any } \gamma <2s\,,\\%[2mm]
\Phi_1(x_0)\, \big(1+\big|\log\Phi_1(x_0)\big|\bigr)  & \qquad\mbox{for }\gamma =2s\,,\\%[2mm]
\Phi_1(x_0)^{2s/\gamma} & \qquad\mbox{for any }\gamma >2s\,.\\
\end{array}\right.
\end{equation}
On the one hand, it is clear  that the upper bounds are sharp when $\gamma<2s$\,, since the powers of $\Phi_1$ in the lower and upper bounds match. On the other hand, when $\gamma \ge 2s$\ and $f$ is just a function in $\LL^\infty$ we cannot expect to have matching powers. Anyway, when dealing with a ``better'' $f$ (for instance if $f$ is zero on $\partial\Omega$), we can significantly improve the above bounds and obtain absolute upper and lower bounds with matching powers also when $\gamma \ge 2s$. We will show that the behaviour of the right hand side at zero is the one that dictates the sharp boundary behavior. This will be a consequence of assumptions (K2) or (K4), depending on the range of parameters.\vspace{-1mm}
\begin{lem}\label{Lem.Green.2}Let $\K$ be the kernel of $\AI$ satisfying $(K2)$, and let $F$ be a positive convex function with $\n(0)=0$.
Let $u:\Omega \to [0,\infty)$ be a bounded measurable function such  that for all $x_0\in\Omega$
\begin{equation}\label{Lem.Green.2.hyp}
\n(u(x_0))\le \kappa_0\int_{\Omega} u(x)\K(x,x_0)\dx.
\end{equation}
Then, there exists a constant $c_{\Omega}>0$ such that for all $x_0\in \Omega$ the following upper bound holds
\begin{equation}\label{Lem.Green.2.est.Upper.II.c}
\n(u(x_0))\le \kappa_0\int_{\Omega} u(x)\K(x,x_0)\dx\le \kappa_1\B_1(\Phi_1(x_0))\,,
\end{equation}
where $\kappa_1=c_{\Omega}\kappa_0\,\n^{-1}\circ\nl(2c_{\Omega}\kappa_0)$, $F^*$ is the Legendre transform of $F$, and $c_{\Omega} $ depends on $s,N,\gamma,\Omega$ only.
\end{lem}
\noindent {\sl Proof of Lemma \ref{Lem.Green.2}.~}We first recall Young inequality for $\n$ and its Legendre transform $\nl$:
\begin{equation}\label{Young.F}
a\,b\, \le \varepsilon \n(a) +\varepsilon \nl\left(\frac{b}{\varepsilon}\right)\qquad\mbox{for all $\varepsilon>0$ and $a,b\ge 0$\,.}
\end{equation}
Combining the above inequality (with $\varepsilon=1/2$) together with \eqref{Lem.Green.2.hyp}\,, we obtain
\begin{equation}\label{Step.0.0.1}\begin{split}
\n(u(x_0))&\le \kappa_0 \int_{\Omega} u(x)\K(x,x_0)\dx
\le \kappa_0  \|u\|_{\LL^\infty(\Omega)} \|\K(x_0,\cdot)\|_{\LL^1(\Omega)}\\
&\le \frac{1}{2}\n\left(\|u\|_{\LL^\infty(\Omega)}\right)
+\frac{1}{2}\nl\bigl(2\kappa_0 c_{4,\Omega}\B_1(\Phi_1(x_0))\bigr),
\end{split}
\end{equation}
where in the last step we have used the Green function estimate \eqref{rem1.1.Green}\,. Taking the supremum over $x_0\in \Omega$ we obtain
\begin{equation}\label{Step.0.0.2}\begin{split}
\n\left(\|u\|_{\LL^\infty(\Omega)}\right)
&\le \nl\left(2 c_{5,\Omega}\kappa_0\right)\qquad\mbox{or}\qquad \|u\|_{\LL^\infty(\Omega)}\le \n^{-1}\circ\nl\left(2 c_{5,\Omega}\kappa_0\right):=\tilde\kappa_1.
\end{split}
\end{equation}
Plugging this last inequality in \eqref{Lem.Green.2.hyp} and using the Green function estimate \eqref{rem1.1.Green}\,, we obtain
\begin{equation*}\begin{split}
\n(u(x_0))&\le \kappa_0 \int_{\Omega} u(x)\K(x,x_0)\dx
\le \kappa_0 \kappa_1 \|\K(x_0,\cdot)\|_{\LL^1(\Omega)}
\le c_{5,\Omega}\kappa_0 \tilde\kappa_1 \B_1(\Phi_1(x_0)),
\end{split}
\end{equation*}
as desired.\qed
\medskip

\noindent\textbf{Remark. } As we shall see later, the bound
\eqref{Lem.Green.2.est.Upper.II.c} is sharp when $\gamma<2s$.

In order to obtain precise estimates near the boundary in the rest of the cases, we need to use the precise behaviour of $\n$ near zero, and this is reflected in the next results.
\begin{lem}\label{Lem.Green.2b}Under the assumptions of Lemma $\ref{Lem.Green.2}$, assume further that for $0\leq a\le 1 $ we have $F(a)\ge \kb a^m$ for some $\kb>0$ and $m>1$, and that for all $x_0\in\Omega$
\begin{equation}\label{Lem.Green.2.hyp.b}
\n(u(x_0))\le \kappa_0\int_{\Omega} u(x)\K(x,x_0)\dx.
\end{equation}
Then, when $\gamma<2sm/(m-1)$\,, there exist a constant $\kappa_4>0$ such that, for all $x_0\in \Omega$,
\begin{equation}\label{Lem.Green.2.est.Upper.II.b1}
\kb u^{m}(x_0) \le \n(u(x_0))\le \kappa_0\int_{\Omega} u(x)\K(x,x_0)\dx\le \kappa_4\kappa_0^{\frac{m}{m-1}}\Phi_1(x_0)\,.
\end{equation}
On the other hand, when $1\ge \gamma\geq  2s m/(m-1)$, there exist a constant $\kappa_5>0$ such that, for all $x_0\in \Omega$,
\begin{equation}\label{Lem.Green.3.est.Lower.II.b}
\begin{split}
\kb u^{m}(x_0) &\le \n(u(x_0))\le \kappa_0\int_{\Omega} u(x)\K(x,x_0)\dx\\
&\le
\kappa_{5} \kappa_0^{\frac{m}{m-1}}
\left\{
\begin{array}{ll}
\Phi_1(x_0)\big(1+|\log \Phi_1(x_0)|^{m/(m-1)}\big)&\text{if $\gamma= 2s m/(m-1)$},\\
\Phi_1(x_0)^{\frac{2sm}{(m-1)\gamma}}&\text{if $\gamma> 2s m/(m-1)$},\\
\end{array}
\right.
\end{split}
\end{equation}
and all the above upper bounds are sharp.
Here, $\kappa_4,\kappa_5>0$ only depend on $s,\gamma, m, N,\Omega$\,.
\end{lem}

\noindent {\sl Proof of Lemma \ref{Lem.Green.2b}.~}In view of Lemma \ref{Lem.Green.2} it is clear that taking $x_0$ close enough to the boundary, we have $u(x_0)\le 1$. Hence, the behaviour of the convex function $\n$ that really matters in the boundary estimates is just the behaviour near zero, namely $\n(u)\ge \kb u^m$\,. The case $\gamma<2s$ then follows by \eqref{Lem.Green.2.est.Upper.II.c}, hence it only remains to  deal with the case $\gamma\ge 2s$

\noindent$\bullet~$\textit{Preliminaries for the case $\gamma\ge 2s$. }The bounds \eqref{Lem.Green.2.est.Upper.II} give the following starting boundary behaviour:
\begin{equation}\label{Step.0.2.Lem.Green.2}\begin{split}
\kb u^m(x_0)&\le \n(u(x_0))
 \le   \kappa_1 \B_1(\Phi_1(x_0))\\
&\le  \kappa_1
\left\{\begin{array}{lll}
\frac{1}{\varepsilon}\Phi_1^{1-\varepsilon}(x_0)\,, & \qquad\mbox{for $\gamma =2s$ and all $\varepsilon\in (0,1]$}\,,\\[2mm]
\Phi_1(x_0)^{\frac{2s}{\gamma}}\,, & \qquad\mbox{for any }\gamma >2s\,.\\
\end{array}\right.
\end{split}
\end{equation}
This behaviour is not sharp and we will improve it through iterations, by splitting different cases. Summing up, at the moment we have
\begin{equation}\label{Step.0.3.Lem.Green.2}
u(x_0)\le \tilde\kappa_1^{1/m}\Phi_1^{\nu_1}(x_0)\,,
\end{equation}
where $\nu_1<1$ satisfies
\begin{equation}\label{Step.0.4.Lem.Green.2}
\nu_1:=\left\{\begin{array}{lll}
%\frac{1}{m}\,, & \qquad\mbox{for any } \gamma <2s\,,\\[2mm]
\frac{1-\varepsilon}{m}\,, & \quad\mbox{for $\gamma =2s$ and all $\varepsilon\in (0,1]$}\,,\\[2mm]
\frac{2s}{m\gamma}\,, & \quad\mbox{for any }\gamma >2s\,.\\
\end{array}\right.
\end{equation}

\noindent$\bullet~$\textit{Iterative step. }We are now going to show the following result, which turns out to be the ``generic'' iterative step that will allow us to prove the results in all the remaining cases, namely when $\gamma\ge 2s$.

\noindent\textsl{Claim. }Assume that for some $0<\nu\leq 1/m$ we have that
\begin{equation}\label{Step.01.hyp.nu}
u(x_0)\le \tilde\kappa\Phi_1^{\nu}(x_0)\qquad\mbox{for all $x_0\in\Omega$}\,.
\end{equation}
Then, when $2s-\gamma(1-\nu)\neq 0$ we have that, for all $x_0\in\Omega$,
\begin{equation}\label{Step.01.claim.nu}
\kb u^{m}(x_0) \le \kappa_0  \int_{\Omega}  u(x)  \K(x,x_0)\dx \leq \bar \kappa_2 \tilde\kappa \kappa_0\biggl(\Phi_1(x_0)^{\nu+\frac{2s}{\gamma}}+ \frac{\Phi_1(x_0) - \Phi_1(x_0)^{\nu+\frac{2s}{\gamma}}}{2s-\gamma(1-\nu)}\biggr),
\end{equation}
where the constant $\bar \kappa_2$ is universal.

\noindent\textit{Proof of the Claim. }Set
$$
\overline{R}:=1+\diam(\Omega)
$$
so that for any $x_0\in\Omega$ we have $\Omega\subseteq B_{\overline{R}}(x_0)$,
fix $x_0\in \Omega$ and choose $R_0$ so that
\[
R_0:=\bar R \Phi_1(x_0)^{1/\gamma}\le \overline{R}.
\]
Recall now the upper part of (K2) assumption, that can be rewritten in the  form
\begin{equation}\label{Step.1.1.Lem.Green.2}
\K(x,x_0)\le
\frac{c_{1,\Omega}}{|x-x_0|^{N-2s}}%\left(\frac{\Phi_1(x)}{|x-x_0|^\gamma}\wedge 1\right)
\left\{\begin{array}{cl}
\dfrac{\Phi_1(x_0)}{|x-x_0|^{\gamma}} &\quad\mbox{for any }x\in \Omega\setminus B_{R_0}(x_0)\\
1&\quad\mbox{for any }x\in B_{R_0}(x_0)\cap \Omega\\
\end{array}\right.
\end{equation}
Next, we recall that $\p=\dist(\cdot,\partial\Omega)^\gamma$\,, and since  $\gamma \leq 1$ we get
\[
|\dist(x,\partial\Omega)-\dist(x_0,\partial\Omega)|\le |x-x_0|\qquad\Longrightarrow\qquad  \p(x)\le\p(x_0)+ |x-x_0|^\gamma.
\]
As a consequence, recalling that $c_0\p(y)\le \Phi_1(y)\le c_1\p(y)$ for all $y\in \overline{\Omega}$, we obtain
\[
\Phi_1(x)\le c_1\p(x)\le c_1\p(x_0)+c_1 |x-x_0|^\gamma\le \frac{c_1}{c_0}\Phi_1(x_0)+c_1 |x-x_0|^\gamma,
\]
so that, for all $x\in\Omega$,
\begin{equation}\label{Step.1.2.Lem.Green.2}
\Phi_1(x)\le \frac{c_1}{c_0}\Phi_1(x_0)+c_1|x-x_0|^\gamma \le \overline{k}\left\{\begin{array}{cl}
|x-x_0|^\gamma  &\quad\mbox{for any }x\in \Omega\setminus B_{R_0}(x_0))\\
\Phi_1(x_0)&\quad\mbox{for any }x\in B_{R_0}(x_0)\cap \Omega\\
\end{array}\right..
\end{equation}
Combining the above estimates, we obtain
\begin{equation}\label{Step.1.3.Lem.Green.2}
\Phi_1(x)^{\nu} \K(x,x_0)\le
\frac{c_{1,\Omega}\overline{k}^{\nu}}{|x-x_0|^{N-2s}}
\left\{\begin{array}{cl}
\dfrac{\Phi_1(x_0)}{|x-x_0|^{\gamma(1-\nu)}} &\quad\mbox{for any }x\in \Omega\setminus B_{R_0}(x_0)\\
\Phi_1(x_0)^{\nu}&\quad\mbox{for any }x\in B_{R_0}(x_0)\cap \Omega\\
\end{array}\right..
\end{equation}
We now recall that $R_0^\gamma= \bar R^\gamma \Phi_1(x_0)$, and that we can assume $\Phi_1(x_0)\le 1/2$ without loss of generality (recall that $x_0$ is a point close to the boundary);   we next use \eqref{Step.1.3.Lem.Green.2} to obtain, for all $2s-\gamma(1-\nu)\neq 0$,
\begin{equation}\label{Step.1.4.Lem.Green.2}\begin{split}
&\int_{\Omega}  u(x)  \K(x,x_0)\dx
\le \tilde\kappa \int_{B_{R_0}(x_0)} \Phi_1^{\nu}(x)\K(x,x_0)\dx+\tilde\kappa \int_{\Omega\setminus B_{R_0}(x_0)} \Phi_1^{\nu}(x)\K(x,x_0)\dx\\
&\le \tilde\kappa c_{1,\Omega}\overline{k}_\gamma^{\nu}\left[\int_{B_{R_0}(x_0)} \frac{\Phi_1^{\nu}(x_0)}{|x-x_0|^{N-2s}} \dx + \int_{B_{\overline{R}}(x_0)\setminus B_{R_0}(x_0)} \dfrac{\Phi_1(x_0)}{|x-x_0|^{N-2s+\gamma(1-\nu)}}\dx\right]\\
&= \tilde\kappa c_{1,\Omega}\overline{k}_\gamma^{\nu}\omega_N\Phi_1^{\nu}(x_0)\left[\frac{R_0^{2s}}{2s}
            + \frac{\Phi_1^{1-\nu}(x_0)}{2s-\gamma(1-\nu)}\left(\overline{R}^{2s-\gamma(1-\nu)}-R_0^{2s-\gamma(1-\nu)}\right)\right]\\
&\leq \bar \kappa_1 \tilde\kappa \biggl(\frac{\Phi_1(x_0)^{\nu+\frac{2s}{\gamma}}}{2s}+ \bar R^{2s-\gamma(1-\nu)}\frac{\Phi_1(x_0) - \Phi_1(x_0)^{\nu+\frac{2s}{\gamma}}}{2s-\gamma(1-\nu)}\biggr)\\
&\leq \bar \kappa_2 \tilde\kappa\biggl(\Phi_1(x_0)^{\nu+\frac{2s}{\gamma}}+ \frac{\Phi_1(x_0) - \Phi_1(x_0)^{\nu+\frac{2s}{\gamma}}}{2s-\gamma(1-\nu)}\biggr),
\end{split}
\end{equation}
Then estimate \eqref{Step.01.claim.nu} follow, using hypothesis \eqref{Lem.Green.2.hyp.b} and recalling that $\n(u)\ge \kb u^m$.
Notice that, in the above estimates, $\bar\kappa_1$ and $\bar\kappa_2$ depend only on $N,s,\gamma,\tilde\kappa, c_{1,\Omega}\overline{k}_\gamma^{\nu}\omega_N,\overline{R}$.

\medskip

Once the iterative step is proven, we need to consider several cases. We shall first consider the case $\gamma>2s$, and finally the case $\gamma=2s$.

\noindent$\bullet~$\textit{Improving the boundary estimates when $\gamma>2s$. }In this case, thanks to \eqref{Step.0.3.Lem.Green.2}, we know that  we can start from a fixed exponent
\[
\nu_1=\frac{2s}{m\gamma}<\frac{1}{m}<1\,,%\qquad\mbox{since}\qquad \frac{2s}{\gamma}<1,
\]
and we want to arrive to the sharp exponent
\[
\nu_\infty=\frac{1}{m}\wedge \frac{2s}{\gamma(m-1)}=\frac{\sigma}{m}
\]
(recall \eqref{eq:sigma}).
We are going to use inequality \eqref{Step.01.claim.nu} with $\nu_1$ as above\,, and we need to split the proof in different cases.

\noindent$\bullet~${\sc Case I. }
We first consider the range
$$
1<\frac{\gamma}{2s} < \frac{m+1}{m}
$$
Observe that, in this range of exponents,
\[
\frac{2s}{\gamma} > \frac{m}{m+1}>\frac{m-1}{m},\qquad\mbox{hence}\qquad \sigma=1\wedge\frac{2sm}{\gamma(m-1)}=1\qquad\mbox{or}\qquad \nu_\infty=\frac{1}{m}\,.
\]
Using the value $\nu_1=2s/m\gamma$ we always have in this case that
\[
2s-\gamma(1-\nu_1)> 0\qquad\mbox{or equivalently}\qquad  \nu_1+\frac{2s}{\gamma}> 1,
\]
hence, recalling \eqref{Step.0.3.Lem.Green.2}, inequality \eqref{Step.01.claim.nu} with $\nu=\nu_1=2s/m\gamma$ gives
\begin{equation}\label{Case.I.01}
\begin{split}
u^{m}(x_0) &\le
\frac{\kappa_0}{\kb}  \int_{\Omega}  u(x)  \K(x,x_0)\dx\\
&  \leq \frac{\bar \kappa_2}{\kb }\kappa_0 \tilde\kappa_1^{1/m} \biggl(\Phi_1(x_0)^{\frac{2s-\gamma(1-\nu_1)}{\gamma}}+ \frac{1}{2s-\gamma(1-\nu_1)}\biggr)\Phi_1(x_0) \le \ka_3\kappa_0\Phi_1(x_0),
\end{split}
\end{equation}
which is the desired upper bound.

\noindent$\bullet~${\sc Case II. }We now proceed with the proof in the range
$$
\frac{\gamma}{2s} > \frac{m+1}{m},\qquad \text{thus}\qquad 2s-\gamma(1-\nu_1)< 0.
$$
We first observe that we can always simplify
\eqref{Step.01.claim.nu} as
\begin{equation}\label{Step.01.claim.nu2}
u^{m}(x_0)\leq \frac{\kappa_0}{\kb}  \int_{\Omega}  u(x)  \K(x,x_0)\dx \le  \bar\kappa_3\tilde \kappa\kappa_0 \frac{ \Phi_1(x_0)^{\nu+\frac{2s}{\gamma}}}{\gamma(1-\nu)-2s}\qquad \text{whenever $2s-\gamma(1-\nu)< 0$}.
\end{equation}
Define the increasing sequence of exponents: (recall that $\nu_1=2s/m\gamma <1/m<1$)
\begin{equation}\label{Case.II.nuk}
\nu_{k}:=\frac{\nu_{k-1}}{m}+\frac{2s}{m\gamma}
=\frac{\nu_1}{m^{k-1}}+\frac{2s}{m\gamma}\left(1+\frac{1}{m}+\frac{1}{m^2}+\dots+\frac{1}{m^{k-2}}\right)\xrightarrow[k\to\, \infty]{}\nu_\infty:=\frac{2s}{\gamma(m-1)}
\end{equation}
Eventually by choosing a smaller $\nu_1$\,, we always have the following alternative:

\noindent$\bullet~$\textsc{Case II.A}) For all $k\ge 1$ we have
\[
2s-\gamma(1-\nu_k)< 0,
\]
in which case
\[
\nu_\infty=\frac{2s}{\gamma(m-1)}\le 1-\frac{2s}{\gamma},\qquad\mbox{hence}\qquad \sigma=1\wedge\frac{2sm}{\gamma(m-1)}=\frac{2sm}{\gamma(m-1)}\,.
\]
\noindent$\bullet~$\textsc{Case II.B}) There exists $\bar k>1$ such that
\[
2s-\gamma(1-\nu_{\bar k})< 0\qquad\mbox{and}\qquad
2s-\gamma(1-\nu_{\bar k+1})> 0.
\]

Note that, starting from
\eqref{Step.01.claim.nu2} with $\nu=\nu_1=2s/m\gamma$ (see \eqref{Step.0.3.Lem.Green.2}), we deduce that
 $$
 u^{m}(x_0) \le  \bar\kappa_3\tilde\kappa_1^{1/m}\kappa_0 \frac{ \Phi_1(x_0)^{\nu_1+\frac{2s}{\gamma}}}{\gamma(1-\nu_1)-2s},
 $$
 hence (see \eqref{Case.II.nuk})
 $$
 u(x_0)\leq \kappa_4\Phi_1(x_0)^{\frac{\nu_1}{m}+\frac{2s}{m\gamma}}=\kappa_4\Phi_1(x_0)^{\nu_2}.
 $$
 In other words, we proved that starting from $\nu_1$ we reach the exponent $\nu_2$ (and more in general the bound \eqref{Step.01.hyp.nu} with $\nu_k$ such that $2s-\gamma(1-\nu_k)<0$ implies that \eqref{Step.01.hyp.nu} holds with $\nu_{k+1}$).
 Hence,
if we are in Case II.B, it means that in a finite number of iteration we fall into Case I, and then the desired result holds by the previous discussion. So we only need to consider Case II.A.

In this case the iteration carries on infinitely many times, so we need to be very careful about how constants enter into the estimates.
We shall distinguish two cases depending whether $\nu_\infty<1-\frac{2s}{\gamma}$ or $\nu_\infty=1-\frac{2s}{\gamma}$.

\noindent
$\bullet~${\sc Case II.A.1)} We assume that $\nu_\infty<1-\frac{2s}{\gamma}$.
In this case the coefficients
$$
A_k:=\frac{1}{\gamma(1-\nu_\k)-2s}
\leq \frac{1}{\gamma(1-\nu_\infty)-2s}=:A_\infty
$$
are uniformly bounded, so starting from \eqref{Step.0.3.Lem.Green.2} and applying
\eqref{Step.01.claim.nu2}
we get
$$
u^m(x_0)\leq
\frac{\kappa_0}{\kb}  \int_{\Omega}  u(x)  \K(x,x_0)\dx\leq A_1\bar\kappa_3\tilde\kappa_0\kappa_1^{1/m}
\Phi_1(x_0)^{\nu_1+\frac{2s}{\gamma}},
$$
hence
$$
u(x_0)\leq \tilde \kappa_2^{1/m}\Phi_1(x_0)^{\nu_2},\qquad \tilde\kappa_2:=A_1\bar\kappa_3\kappa_0\tilde\kappa_1^{1/m}.
$$
More in general, iterating $k$ times we get
$$
u(x_0)^m\leq \frac{\kappa_0}{\kb}  \int_{\Omega}  u(x)  \K(x,x_0)\dx\leq \tilde \kappa_k\Phi_1(x_0)^{m\nu_k},\qquad \tilde \kappa_k:=A_{k-1}\bar\kappa_3\kappa_0\tilde\kappa_{k-1}^{1/m}.
$$
Note that as $k\to \infty$ the sequence $\tilde\kappa_k$ remains bounded and converges to $\tilde\kappa_\infty\approx (A_\infty\bar\kappa_3\kappa_0)^{m/(m-1)}$,
concluding the proof of this case.

\noindent
$\bullet~${\sc Case II.A.2)} We assume that $\nu_\infty=1-\frac{2s}{\gamma}=\frac{1}{m}$.
In this case $\nu_k\to \nu_\infty= 1/m$ and $A_k:=\frac{1}{\gamma(1-\nu_\k)-2s}\to \infty$ as $k\to \infty.$
Thus, fixed a point $x_0$, we argue as in Case II.A.1
for $k \leq k_0$, where $k_0$ will be chosen below.
In this way we get
$$
u(x_0)^m\leq \frac{\kappa_0}{\kb}  \int_{\Omega}  u(x)  \K(x,x_0)\dx\leq \tilde \kappa_{k_0}\Phi_1(x_0)^{m\nu_{k_0}},\qquad \tilde \kappa_{k}:=A_{k-1}\bar\kappa_3\kappa_0\tilde\kappa_{k-1}^{1/m}.
$$
Noticing that $A_k=\frac{1}{\gamma(\nu_\infty-\nu_k)}\leq \frac{1}{\gamma(\nu_\infty-\nu_{k_0})}$
we can bound
$$
\tilde \kappa_{k_0}\leq \bar \kappa_\infty A_{\kappa_0-1}A_{k_0-2}^{1/m} \cdots A_{1}^{1/m^{k_0}}\leq \hat \kappa_\infty (\nu_\infty-\nu_{k_0})^{-m/(m-1)},
$$
therefore (recall that in this case $\nu_\infty=1/m$)
$$
u(x_0)^m\leq \frac{\kappa_0}{\kb}  \int_{\Omega}  u(x)  \K(x,x_0)\dx\leq
\hat \kappa_\infty\Phi_1(x_0) \frac{\Phi_1(x_0)^{-m(\nu_\infty-\nu_{k_0})}}{(\nu_\infty-\nu_{k_0})^{m/(m-1)}}.
$$
Now, if we choose $k_0$ large enough so that $\nu_\infty-\nu_{k_0}\approx |\log \Phi_1(x_0)|^{-1}$ then
$$
\Phi_1(x_0)^{-m(\nu_\infty-\nu_{k_0})}=
e^{-m(\nu_\infty-\nu_{k_0})\log \Phi_1(x_0)}\approx 1
$$
and we get
$$
u(x_0)^m\leq \frac{\kappa_0}{\kb}  \int_{\Omega}  u(x)  \K(x,x_0)\dx\leq
\bar \kappa_\infty\Phi_1(x_0) |\log \Phi_1(x_0)|^{m/(m-1)},
$$
as desired.

\noindent$\bullet~${\sc Case III. }
Assume that $\gamma/2s= (m+1)/m$. In this case
it suffices to slightly reduce the value of $\nu_1$ to fall into Case II, and because in this case $\nu_\infty+\frac{2s}{\gamma}>1$, we fall into Case II.B.

\noindent$\bullet~${\sc Case IV. }We now assume that $\gamma=2s$.
In this case we have $\nu_1=(1-\varepsilon)/m$ for any $\varepsilon\in (0,1)$, and we can always choose $\varepsilon>0$ small enough so that we fall into Case I.\qed

\medskip

We can now prove Proposition \ref{cor.AI.Phi1}.

\noindent {\sl Proof of Proposition \ref{cor.AI.Phi1}.~}The proof begins by noticing that, using the kernel representation of $\AI$, we have
\begin{equation}\label{Prop72.1}
\Phi_1(x)=\lambda_1\,\AI\Phi_1(x)=\lambda_1\int_\Omega \Phi_1(y)\K(x,y)\dy.
\end{equation}
\noindent$\bullet~$\textsc{Upper bounds. }We first derive the sharp upper bounds. The above formula suggests to repeat the proof of Lemma \ref{Lem.Green.2b} above with $m=1$.

 Since we know that $\Phi_1\in \LL^\infty(\Omega)$ by Proposition \ref{prop.AI.Phi1}\,, as a consequence of \eqref{Prop72.1} we obtain
\begin{equation}\label{Prop72.2}
\Phi_1(x)=\lambda_1\int_\Omega \Phi_1(y)\K(x,y)\dy\le c_{4,\Omega}\lambda_1\|\Phi_1\|_{\LL^\infty(\Omega)}\B_1(\p(x)) \le \kappa_1 \B_1(\p(x)),
\end{equation}
where we used the estimates on the $\LL^1$ norm of $\K$ given by Lemma \ref{Lem.Green}.

Recalling \eqref{rem1.2.Green}
we see that the upper bounds are already sharp when $\gamma<2s$, so we need to prove the case $\gamma\ge 2s$.
We begin by repeating the first steps of the proof of Lemma \ref{Lem.Green.2b}, namely we  follow the proof of Lemma \ref{Lem.Green.2b}, letting $m=1$ there, and we obtain the analogous of formula \eqref{Step.0.2.Lem.Green.2}, that we rewrite in the form analogous to formula \eqref{Step.0.3.Lem.Green.2} as follows: for all $x_0\in \Omega$\,, we have
\begin{equation}\label{Prop72.4}
\Phi_1(x_0)\le  \kappa_1 \delta^{\gamma\nu_1}(x_0)\,,
\end{equation}
with $\nu_1<1$ given by
\begin{equation}\label{Prop72.5}
\nu_1:=\left\{\begin{array}{lll}
1-\varepsilon\,, & \quad\mbox{for $\gamma =2s$ and all $\varepsilon\in (0,1]$}\,,\\[2mm]
\frac{2s}{\gamma}\,, & \quad\mbox{for any }\gamma >2s\,.\\
\end{array}\right.
\end{equation}
This behaviour is not sharp
but we can improve it exactly as in the proof of Lemma \ref{Lem.Green.2b}
(see
the proof of \eqref{Step.01.claim.nu} starting from
\eqref{Step.01.hyp.nu})
noticing that in this case, since $m=1$,
at each iteration the exponent improves by a fixed amount (equivalently $\nu_\infty=\infty$, see \eqref{Case.II.nuk}),
so (up to slightly reducing the value of $\nu_1$) we always fall into Case I or Case II.B.

\noindent$\bullet~$\textsc{Lower bounds. }The upper bound proven above, namely $\Phi_1\le \ka \p$ on $\Omega$\,, together with formula \eqref{Prop72.1} and assumption (K2), implies the desired sharp lower bound: more precisely, the lower bound in assumption (K2) gives, for all $x\in \Omega$,
\[
\Phi_1(x)=\lambda_1(\AI\Phi_1)(x)\ge c_{0,\Omega}\,\lambda_1\p(x)\int_\Omega \Phi_1(y)\p(y)\dy \ge \frac{c_{0,\Omega}\lambda_1}{\ka}\,\p(x)\int_\Omega \Phi_1(y)^2\dy
=\kb \,\p(x)
\]
where $\kb=c_{0,\Omega}\lambda_1/\ka$ and we used that $\|\Phi_1\|_{\LL^2(\Omega)}=1$.\qed

%%%%%%%%%%%%%%%%%%%%%%%%%%%%%%%%%%%%%%%%%%%%%%%%%%%%%%%%%%%%%%%%%%%

\subsubsection{Green function estimates II. The lower bounds}
In the same spirit as in the previous section, we now obtain sharp lower estimates.

\begin{lem}\label{Lem.Green.3}Let $\K$ be the kernel of $\AI$ satisfying $(K2)$, and let $F$ be a nonnegative increasing function with $\n(0)=0$\,.
Let $u:\Omega \to [0,\infty)$ be a bounded measurable function, assume that  $F(r)\le \ka r^m$ for some $m>1$ for all $0\leq r\leq 1$, and that for all $x_0\in\Omega$
\begin{equation}\label{Lem.Green.3.hyp}
\n(u(x_0))\ge \kappa_0\int_{\Omega} u(x)\K(x,x_0)\dx.
\end{equation}
Then, there exists a constant $\kappa_{6}>0$ such that, for all $x_0\in \Omega$,
\begin{equation}\label{Lem.Green.3.est.Lower.II}
  \ka u^m(x_0) \ge \n(u(x_0))
\ge \kappa_0\int_{\Omega} u(x)\K(x,x_0)\dx\ge \kappa_6\kappa_0^{\frac{m}{m-1}} \Phi_1(x_0)\,,
\end{equation}
and the estimates are sharp when $\gamma< 2s m/(m-1)$\,.

On the other hand, when $1\ge \gamma\geq  2s m/(m-1)$, (K4) implies that there exist a constant $\kappa_7>0$ such that, for all $x_0\in \Omega$,
\begin{equation*}
\begin{split}
  \ka u^m(x_0) &\ge \n(u(x_0))
\ge \kappa_0\int_{\Omega} u(x)\K(x,x_0)\dx\\
&\ge
\kappa_{7} \kappa_0^{\frac{m}{m-1}}
\left\{
\begin{array}{ll}
\Phi_1(x_0)\big(1+|\log \Phi_1(x_0)|^{m/(m-1)}\big)&\text{if $\gamma= 2s m/(m-1)$},\\
\Phi_1(x_0)^{\frac{2sm}{(m-1)\gamma}}&\text{if $\gamma> 2s m/(m-1)$},\\
\end{array}
\right.
\end{split}
\end{equation*}
and the estimates are sharp. Here, $\kappa_6,\kappa_7>0$ only depend on $s,\gamma, m, N,\Omega$\,.
\end{lem}
\noindent {\sl Proof.~}As already discussed, since we are interested in boundary behaviour, we can assume without loss of generality that $u(x_0)\le 1$.

Recalling that $\ka u^m(x_0) \ge \n(u(x_0))$ when $0\le u(x_0)\le 1$, it follows from (K2) that
\begin{equation}\label{Lem.Green.3.1}\begin{split}
\ka u^{m}(x_0) \ge \n(u(x_0))&\ge \kappa_0\int_\Omega u(x)\K(x,x_0)\dx\\
&\ge \kappa_0 c_{0,\Omega} \Phi_1(x_0)\int_\Omega u(x)\Phi_1(x)\dx = \kappa_0c_{0,\Omega} \Phi_1(x_0) \|u\|_{\LL^1_{\Phi_1}(\Omega)}\,.
\end{split}
\end{equation}
Multiplying the above lower bound by $\Phi_1(x_0)$ and integrating with respect to $x_0\in \Omega$, we obtain an absolute bound for the $\LL^1$ weighted norm:
\begin{equation}\label{Lem.Green.3.2}\begin{split}
\|u\|_{\LL^1_{\Phi_1}(\Omega)}^{\frac{m-1}{m}} \ge \left(\frac{\kappa_0}{\ka}c_{0,\Omega}\right)^{\frac{1}{m}} \|\Phi_1\|_{\LL^{\frac{m+1}{m}}(\Omega)}^{\frac{m+1}{m}}\,.
\end{split}
\end{equation}
We now reinsert the above absolute lower bound \eqref{Lem.Green.3.2} into \eqref{Lem.Green.3.1} to obtain
\begin{equation}\label{Lem.Green.3.3}\begin{split}
\ka u^m(x_0) \ge \n(u(x_0))&\ge \kappa_0\int_\Omega u(x)\K(x,x_0)\dx\\
&\ge \kappa_0c_{0,\Omega} \Phi_1(x_0) \|u\|_{\LL^1_{\Phi_1}(\Omega)}
\ge c_{6,\Omega}\kappa_0^{\frac{m}{m-1}} \Phi_1(x_0)\,,
\end{split}
\end{equation}
where $c_{6,\Omega}\sim (c_{0,\Omega}/\ka)^{1/(m-1)}$\,; this is exactly \eqref{Lem.Green.3.est.Lower.II}\,, and this lower bound is sharp  when $\gamma< 2s m/(m-1)$\,, in view of the matching upper bounds \eqref{Lem.Green.2.est.Upper.II} of Lemma \ref{Lem.Green.2}.

It only remains to prove the sharp lower bounds when $\gamma\geq 2s m/(m-1)$ and for this reason from now on we will assume (K4): for all $x,y\in\overline{\Omega}$
\begin{equation}\label{Lem.Green.3.K4}
\K(x,y)\ge \frac{c_{0,\Omega}}{|x-y|^{N-2s}}
\left(\frac{\Phi_1(x)}{|x-y|^\gamma}\wedge 1\right)
\left(\frac{\Phi_1(y)}{|x-y|^\gamma}\wedge 1\right)\,.
\end{equation}
As a consequence, letting $R_0^\gamma=\p(x_0)\asymp\Phi_1(x_0)$ we get that (changing $c_{0,\Omega}$ if needed)
\begin{equation}\label{Lem.Green.3.K4.1}
\K(x,x_0)\ge \frac{c_{0,\Omega}}{|x-x_0|^{N-2s}}\,,\qquad\mbox{for all }x\in B_{R_0/2}(x_0)\,.
\end{equation}
Note that, by the estimate proved above, if we set $\nu_1=1/m$ and $\kappa_1=c_{6,\Omega}^{1/m}\kappa_0^{\frac{1}{m-1}}$, then we have
\[
u(x)\ge \kappa_1\Phi_1(x)^{\nu_1}\ge \kappa_1\tilde{\kappa}_\gamma\Phi_1(x_0)^{\nu_1}\,,\qquad\mbox{for all }x\in B_{R_0/2}(x_0)
\]
(note that $|x-x_0|<R_0/2$ and $\Phi_1(x)\ge \kappa_\gamma\p(x)\ge \kappa_\gamma(\p(x_0)-|x-x_0|^\gamma)\ge \kappa_\gamma\p(x_0)/2\ge \tilde{\kappa}_\gamma\Phi_1(x_0)$).

Next we observe that, proceeding analogously to \eqref{Step.1.4.Lem.Green.2}, for all $2s-\gamma(1-\nu_1)\neq 0$ we obtain
\begin{equation}\label{Step.1.4.Lem.Green.3}\begin{split}
\int_{\Omega}  u(x)  \K(x,x_0)\dx
&\ge \kappa_1^{\nu_1} \tilde{\kappa}_\gamma^{\nu_1}\Phi_1(x_0)^{\nu_1} \int_{B_{R_0}(x_0)} \K(x,x_0)\dx\\
&\ge \kappa_1^{\nu_1} \tilde{\kappa}_\gamma^{\nu_1}\Phi_1(x_0)^{\nu_1} \int_{B_{R_0/2}(x_0)} \frac{1}{|x-x_0|^{N-2s}} \dx\\
&= \kappa_1^{\nu_1} \frac{\tilde{\kappa}_\gamma^{\nu_1}}{2s}\Phi_1(x_0)^{\nu_1} R_0^{2s}
= \kappa_1^{\nu_1} \frac{\tilde{\kappa}_\gamma^{\nu_1}}{2s}\kappa_1^{\frac{1}{\nu_1}}\Phi_1(x_0)^{\nu_1+\frac{2s}{\gamma}}\,.
\end{split}
\end{equation}
Thus, since $\ka u^{m}(x_0) \ge \n(u(x_0))$ when $0\le u(x_0)\le 1$, it follows
\begin{equation}\label{Lem.Green.3.3b}\begin{split}
u (x_0) \ge \kappa_2 \kappa_1^{\frac{1}{m\nu_1}}\Phi_1(x_0)^{\frac{\nu_1}{m}+\frac{2s}{\gamma m}}:=\kappa_2 \kappa_1^{\frac{1}{m\nu_1}}\Phi_1(x_0)^{\nu_2}\,.
\end{split}
\end{equation}
Iterating this process, analogously to Case II.A of Lemma \ref{Lem.Green.2b}, if we set
\[
\nu_{k}:=\frac{\nu_{k-1}}{m}+\frac{2s}{m\gamma}
=\frac{\nu_1}{m^{k-1}}+\frac{2s}{m\gamma}\left(1+\frac{1}{m}+\frac{1}{m^2}+\dots+\frac{1}{m^{k-2}}\right)
\]
the iterative step becomes
\begin{equation}\label{Step.1.9.Lem.Green.2b}
u(x_0)\ge \tilde{\kappa}_{k}\kappa_0^{\frac{1}{m-1}}\Phi_1^{\nu_k}(x_0)\qquad\mbox{implies}\qquad u(x_0)\ge \tilde{\kappa}_{k+1}\kappa_0^{\frac{1}{m-1}}\Phi_1^{\nu_{k+1}}(x_0)\,.
\end{equation}
Letting $k\to \infty$ we get $\nu_k\to \nu_\infty=2s/[\gamma(m-1)]$ and
\begin{equation}\label{Step.1.10.Lem.Green.2b}
u(x_0)\ge \kappa_\infty\, \kappa_0^{\frac{1}{m-1}}\Phi_1^{\nu_\infty}(x_0)\,,
\end{equation}
where the constant $\kappa_\infty\sim \tilde{\kappa}_2^{1/(m-1)}$ is positive and depends only on $m,s,\gamma,N,\Omega$. This concludes the proof also in the case $\gamma>2sm/(m-1)$.

We now deal with the ``critical'' case  $\gamma=2sm/(m-1)$, where in view of the upper bound we need to improve our estimates by a logarithmic factor.

\begin{figure}[h]
\center \includegraphics[scale=0.3]{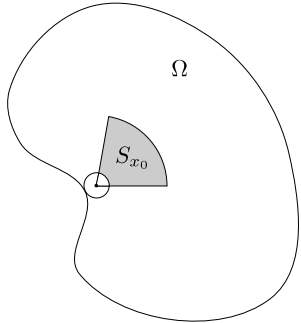}
\caption{The set $S_{x_0}$.}
\label{fig:sector}
\end{figure}

To this aim, we fix $x_0$ close to $\partial\Omega$, we set $R_0:=\Phi_1(x_0)^{1/\gamma}$, and we define the following set
$$
S_{x_0}:=\left\{x \in \Omega\cap B_\rho(x_0)\,:\,|x-x_0|\geq R_0,\, (x-x_0)\cdot \nu_{x_0}\geq \textstyle{\frac23}|x-x_0|\right\},
$$
where $\nu_{x_0}\in \mathbb S^{n-1}$ is a vector pointing in the interior of $\Omega$ and $\rho$ is a small number (still, much larger that $R_0$) depending only on the geometry of $\Omega$ (but independent of $x_0$, provided $x_0$ is close enough to the boundary) to ensure that
\begin{equation}
\label{eq:sector}
S_{x_0}\subset \{x \in \Omega\,:\,\Phi_1(x_0)\leq |x-x_0|^\gamma \leq \hat \kappa\Phi_1(x)\},
\end{equation}
see Figure \ref{fig:sector} (recall that $\Phi_1\approx \delta^\gamma$,
and note that $|x-x_0|^\gamma\leq \delta^\gamma(x)$ inside $S_{x_0}$).
We now start from the bound established before, namely $u \geq \tilde\kappa_1 \Phi_1^{1/m}$ (see \eqref{Lem.Green.3.3}).
Note that thanks to (K4) and \eqref{eq:sector} we have
\begin{equation}\label{Lem.Green.3.K4.sector}
\K(x,x_0)\ge \frac{c_{0,\Omega} \Phi_1(x_0)}{|x-x_0|^{N-2s+\gamma}}\,,\qquad\mbox{for all }x\in S_{x_0}\,.
\end{equation}
Also, again by \eqref{eq:sector}
we deduce that
\begin{equation}\label{Lem.u.sector}
u(x) \geq \tilde\kappa_1 \Phi_1(x)^{1/m}\geq \frac{\tilde\kappa_1}{\hat \kappa^{1/m}}|x-x_0|^{\gamma/m} \,,\qquad\mbox{for all }x\in S_{x_0}\,.
\end{equation}
Hence, applying \eqref{Lem.Green.3.hyp}, \eqref{Lem.Green.3.K4.sector},
and \eqref{Lem.u.sector}, we obtain (recall that $R_0 \ll \rho$)
\begin{align*}
\ka u^m(x_0) \ge \n(u(x_0)) &\geq \kappa_0 \int_{S_{x_0}}u(x)\K(x,x_0)\dx
\geq \kappa_0c_{0,\Omega}\frac{\tilde\kappa_1}{\hat \kappa^{1/m}}\Phi_1(x_0)\int_{S_{x_0}}\frac{1}{|x-x_0|^{N-2s+\gamma+\gamma/m}}\dx\\
&=\kappa_0c_{0,\Omega}\frac{\tilde\kappa_1}{\hat \kappa^{1/m}}\Phi_1(x_0)\int_{S_{x_0}}\frac{1}{|x-x_0|^{N}}\dx=c_N\kappa_0c_{0,\Omega}\frac{\tilde\kappa_1}{\hat \kappa^{1/m}}\Phi_1(x_0)\int_{R_0}^\rho \frac{1}{t}\dt\\
&\geq \kb\kappa_0 \tilde\kappa_1 \Phi_1(x_0)|\log R_0| =\kb \kappa_0\tilde\kappa_1   \gamma \Phi_1(x_0)|\log \Phi_1(x_0)|
\end{align*}
for some universal constant $\kb>0$,
which proves that $u \geq (\kb \kappa_0\tilde\kappa_1\gamma/\ka)^{1/m} \Phi_1^{1/m}|\log \Phi_1|^{1/m}$.

To gain the optimal exponent on the logarithm we need to iterate this estimate:
more precisely, once we know that $u \geq \tilde\kappa_1 \Phi_1^{1/m}|\log \Phi_1|^{\alpha}$ for some $\alpha<\frac{1}{m-1}$, then
instead of \eqref{Lem.u.sector} we have
$$
u(x) \geq \tilde\kappa_1 \Phi_1(x)^{1/m}|\log \Phi_1(x)|^{\alpha}\geq \kb_0 \tilde\kappa_1 |x-x_0|^{\gamma/m}|\log |x-x_0||^{\alpha} \,,\qquad\mbox{for all }x\in S_{x_0}\,,
$$
thus
\begin{align*}
\ka u^m(x_0)  &\geq \kappa_0 \int_{S_{x_0}}u(x)\K(x,x_0)\dx
\geq \kappa_0\kb_1\tilde \kappa_1\Phi_1(x_0)\int_{S_{x_0}}\frac{|\log |x-x_0||^{\alpha}}{|x-x_0|^{N}}\dx\\
&\geq \kb_2\kappa_0 \tilde \kappa_1 \Phi_1(x_0)|\log R_0|^{1+\alpha} = \kb_2\kappa_0 \tilde \kappa_1\gamma \Phi_1(x_0)|\log \Phi_1(x_0)|^{1+\alpha},
\end{align*}
that is $u \geq (\kb_2\kappa_0 \tilde\kappa_1\gamma/\ka)^{1/m}\Phi_1^{1/m}|\log \Phi_1|^{(1+\alpha)/m}$ for some universal constant $\kb_2$.
Iterating infinitely many times, as before the constants do not blow up and behave as $\kappa_0^{\frac{m}{m-1}}$, while the exponents in the logarithm converge to $1/(m-1)$, as desired.\qed

\subsection{Proof of Theorems \ref{Thm.upper.ell} and \ref{Thm.Elliptic.Harnack.m}.~}
\label{section:proof thms}
The proof is now easy in view of the previous results, and can be split into two steps.

\noindent$\bullet~$\textit{A pointwise equality. } Set $V:=f(u)$. We claim that, for a.e. $x_0\in\Omega$, we have
\begin{equation}\label{Step.1.1Thm.Elliptic.Harnack.ineq}
\n(V(x_0))=\int_{\Omega}V(x) \,\K(x_0,x)\,\dx\,.
\end{equation}
To prove this formula, we first use Definition \ref{Def.Very.Weak.Sol.Dual.elliptic} to get
$$
\int_{\Omega} \n(V)\,\psi\,\dx=\int_{\Omega}\AI (V) \,\psi\,\dx=\int_{\Omega}V \,\AI\psi\,\dx
$$
for any $\psi$ such that $\psi/\Phi_1\in \LL^\infty(\Omega)$\,. The proof of formula \eqref{Step.1.1Thm.Elliptic.Harnack.ineq} now follows by considering a sequence of admissible test functions $\psi_{n}^{(x_0)}$ that converge to $\delta_{x_0}$, cp. Step 4 in the proof of Theorem 5.1 of \cite{BV-PPR2-1}\,.

\noindent$\bullet~$\textit{Application of the Green function estimates. }
Thanks to \eqref{Step.1.1Thm.Elliptic.Harnack.ineq}, we apply
\eqref{Lem.Green.2.est.Upper.aaa} to $V=f(u)$
to prove Theorem \ref{Thm.upper.ell}.

For Theorem \ref{Thm.Elliptic.Harnack.m}, we apply Proposition \ref{Prop.Green.2aaa} (with $\kappa_0=1$) to the function $V=f(u)$; in this way, when $\gamma<2sm/(m-1)$ we obtain, for a.e. $x_0\in\Omega$,
\[%\begin{equation}\label{Step.2.1Thm.Elliptic.Harnack.ineq}
 \kappa_6 \Phi_1(x_0)\le V^m(x_0)\le  \kappa_4 \Phi_1(x_0)\,.
\]%\end{equation}
Analogously, when $1\ge \gamma> 2s m/(m-1)$\,, assuming (K4)\,,  we have
\[
\kappa_{7}\Phi_1(x_0)^{\frac{2sm}{(m-1)\gamma}}\le V^m(x_0)\le\kappa_{5}\Phi_1(x_0)\Phi_1(x_0)^{\frac{2sm}{(m-1)\gamma}}\,,\qquad\mbox{for a.e. $x_0\in \Omega$.}
\]
Finally, when $1\ge \gamma= 2s m/(m-1)$\,, assuming (K4)  we have
\[
\kappa_{7}\Phi_1(x_0)\big(1+|\log \Phi_1(x_0)|^{\frac{m}{m-1}}\big)\le V^m(x_0)\le\kappa_{5}\Phi_1(x_0)\big(1+|\log \Phi_1(x_0)|^{\frac{m}{m-1}}\big)\,,\quad\mbox{for a.e. $x_0\in \Omega$.}
\]
The result follows by noticing that $V^m\asymp F(V)=u$ (recall that $m=1/p$).
We finally recall that the constants $\kappa_4,\dots,\kappa_7>0$ depend on $s,\gamma, m, N, \kb,\ka, \Omega$ as in Proposition \ref{Prop.Green.2aaa}\,. \qed

\section{Regularity}\label{sec.regularity}
In this section we discuss the regularity of solutions, both in the interior and up to the boundary. Some proofs follow the same lines of the parabolic results of the companion paper \cite{BFV-Parabolic}, but since here we are dealing with a more general nonlinearity we sketch them below.

\noindent In order to obtain the regularity results, we assume that the operator $\A$ satisfies (A1) and (A2), and we need some assumption on the kernel $K$ of the operator (rather than on its inverse). Such assumptions are quite general, and essentially all the examples at hand do satisfy them.
\subsection{Interior Regularity}
We prove the following result:
\begin{thm}\label{thm.regularity.1}Assume that $\A$ satisfies  $(A1)$ and $(A2)$\,, and that the operator has the following representation formula:
\begin{equation*}%\label{Operator.Hyp.Calpha}
\A g(x)=P.V.\int_{\RR^N}\big(g(x)-g(y)\big)K(x,y)\dy+ D(x)g(x)\,,
\end{equation*}
with
$$
K(x,y)\asymp |x-y|^{-(N+2s)}\quad \text{ for $x,y\in B_r(x_0)$, \, \, where $B_{2r}(x_0)\subset \Omega$},
$$
$$
K(x,y)\lesssim |x-y|^{-(N+2s)}\quad \text{ otherwise},
$$
and $|D(x)|\lesssim \dist(x,\partial\Omega)^{-2s}$.
Let $u$ be a nonnegative bounded weak dual solution to problem \eqref{Elliptic.prob} such that
\begin{align*}
0\leq u(x)\leq M\qquad\mbox{for a.e. $x\in \Omega$, for some $M>0$.}
\end{align*}
\begin{enumerate}[leftmargin=*]
\item[(i)] Let $K(x,y)=K(y,x)$. Then $u$ is {H\"older continuous in the interior}. More precisely, there exists $\alpha>0$ such that,
\begin{equation}\label{thm.regularity.1.bounds.0}
    \|u\|_{C^{\alpha}(B_{r}(x_0))}\leq C.
\end{equation}
\item[(ii)]
Set $\tilde K(x,z):=K(x,z+x)$, and
assume that $\tilde K(x,z)=\tilde K(x,-z)$. Suppose
in addition that
$|\tilde K(x,z)-\tilde K(x',z)|\le c |x-x'|^\beta\,|z|^{-(N+2s)}$ for some $\beta\in (0,1\wedge 2s)$  is not an integer,
that $|D^2_z\tilde K(x,z)|\le C\,|z|^{-(N+2+2s)}$, that $D\in C^\beta(B_{3r/2}(x_0))$, and that $f\in C^\beta(\RR)$.  Then $u$ is a classical solution in the interior:
\begin{equation}\label{thm.regularity.1.bounds.1}
\|u\|_{C^{2s+\beta'}(B_{r}(x_0))}\le C,
\end{equation}
\end{enumerate}
where $\beta'>0$ is any exponent less than $2s\frac{\beta}{1-\beta}$.
Here, the constants in the regularity estimates depend on the solution only through $M>0$.
\end{thm}

\noindent\textbf{Remarks. }(i) In order to guarantee that solutions are  bounded, we need to ensure that Theorem \ref{Thm.upper.ell} holds.
Thus, to apply the regularity results above, one needs to assume that $\AI$ satisfies (K2)
and that $\n(u)=f^{-1}(u)$ is a positive convex function with $F(0)=0$.

\noindent(ii) The above theorem applies to all the examples mentioned in this paper: while this holds for the Restricted and Censored Fractional Laplacian, which have an explicit kernel, it is not so obvious for the Specral Fractional Laplacian. Still, it has been proven in \cite{SV2003} that, when $\A$ is the Spectral Fractional Laplacian, then it can be expressed in the form
\[
\A g(x)=P.V.\int_{\RR^N} \big(g(x)-g(y)\big)\,K(x,y)\dy + D(x)g(x)
\]
with a symmetric kernel $K(x,y)$ supported in $\overline{\Omega}\times\overline{\Omega}$ and satisfying
\[
K(x,y)\asymp\frac{1}{|x-y|^{N+2s}}
\left(\frac{\Phi_1(x)}{|x-y|^\gamma }\wedge 1\right)
\left(\frac{\Phi_1(y)}{|x-y|^\gamma }\wedge 1\right)\quad\mbox{and}\quad D(x)\asymp \Phi_1(x)^{-\frac{2s}{\gamma}}\,.
\]
 We note that interior and boundary regularity estimates for linear and semilinear Dirichlet problems for SFL-type operators have been also obtained in \cite{CS2016,CDDS} using different methods.

\medskip

\noindent\textit{Proof of Theorem \ref{thm.regularity.1}. }
The strategy to prove Theorem \ref{thm.regularity.1} relies on elliptic regularity for linear nonlocal equations. More precisely,  interior H\"older regularity will follow by applying $C^{\alpha}$ estimates of \cite{Ka}  for a ``localized'' linear problem. Similar estimates, valid for viscosity solutions to fully nonlinear equations that may apply to our case, are contained in \cite{BCI, CS2, CS3,  Si}. Once H\"older regularity is established, under a mild H\"older continuity assumption on the kernel (and on $f$) we can use the   Schauder estimates proved in  \cite{BaFiVa} to conclude. Similar Schauder estimates have been proven also in \cite{DK}.
\medskip

\noindent$\bullet$~\textit{Localization of the problem. }
Up to a rescaling, we can assume $r=2$. Take a cutoff function $\rho\in C^\infty_c(B_2)$ such that $\rho\equiv 1$ on $B_1$, and define $v=\rho u$. By construction $u=v$ on $B_1$, so that we can write the equation for $v$ on the small ball $B_1$ as
\begin{equation}\label{localization}
f(v(x)) = \A [v](x) +g_0(x),
\end{equation}
where
$$
g_0(x):= \A\left[(1-\rho)u\right](x)+D(x)u(x)=\int_{\RR^N\setminus B_1}(1-\rho(y))u(y)  K(x,y)\dy+D(x)u(x).
$$
(Note that $(1-\rho)u\equiv 0$ in $B_1$.)

\noindent$\bullet$~\textit{H\"older continuity in the interior. }%\label{sssect.Calpha}
Set $b:= f(v)-g_0$, with $g_0$ as above. Equation \eqref{localization} reads then $\A[v](x)=b(x)$, with $b\in \LL^\infty(B_1)$. Indeed,  it is easy to check that $g_0 \in \LL^\infty(B_1)$, since $K(x,y)\lesssim |x-y|^{-(N+2s)}$ and $D$ is bounded inside $B_1$; moreover, by our assumptions on $f$ and since $0\le u\le M$ on $B_2$\,, we have $f(v)=f(u)\in \LL^\infty(B_1)$. Recalling that $\|v\|_{\LL^\infty(\RR^N)}\le \|u\|_{\LL^\infty(B_2)}<+\infty$, we are now in the position to apply the H\"older estimate of Theorem 1.1 of \cite{Ka} (or also the results of \cite{BCI, CS2, CS3, Si}) to ensure that
\[
\|v\|_{C^{\alpha}(B_{1/2})}\leq C\bigl(\|b\|_{\LL^\infty(B_1)}+\|v\|_{\LL^\infty(\RR^N)}\bigr)
\]
for some universal exponent $\alpha>0$. This proves Theorem \ref{thm.regularity.1}(i).

\medskip

\noindent$\bullet$~\textit{Classical solutions in the interior. }
Under the assumptions on $K$ in part (ii) of the theorem, we can use \cite{JJ} to obtain that the H\"older regularity of $u$ still holds.
Then, once we know that $u\in C^{\alpha}(B_{1/2})$, we repeat the localization argument above with the cutoff function $\rho$  supported inside $B_{1/2}$ to ensure that $v:=\rho u$ is H\"older continuous in $\RR^N$, namely $\|v\|_{{C^{\alpha}(\RR^N)}}<+\infty$. Then, to obtain higher regularity we argue as follows.

Set $\beta_1:=\alpha\beta$. By our assumptions on $K$ and $D$ we have  $g_0\in C^{\beta_1}(B_{1/2})$. Thus, by the previous part of the theorem and since $f\in C^\beta(\RR)$ we can conclude that $b=g_0+f(u)\in C^{\beta_1}(B_{1/2})$. This allows us to apply the Schauder estimates of Theorem 5 of \cite{BaFiVa} (see also \cite{DK,JJ})  to obtain that
\begin{equation}\label{schauder.f}
\|v\|_{C^{2s+\beta_1}(B_{1/4})}
\leq C\bigl(\|b\|_{C^{\beta_1}(B_{1/2})}+\|v\|_{{C^{\beta_1}(\RR^N)}}\bigr).
\end{equation}
 In case $\beta_1=\beta$ we stop here. Otherwise we set
$\alpha_1:=\min\{1,2s+\beta\}$ and we repeat the argument above with $\beta_2:=\alpha_1\beta$
in place of $\beta_1$. In this way, we obtain that $v\in C^{2s+\beta_1}(B_{2^{-4}})$. Iterating this procedure finitely many times, one can reach any exponent $\beta'$ smaller than $2s\frac{\beta}{1-\beta}$.
Finally, a covering argument completes the proof of Theorem \ref{thm.regularity.1}(ii).\qed

\subsection{Boundary regularity}

\medskip

We now prove H\"older regularity up to the boundary under the assumption that $f\in C^\beta(\RR)$ is a nonnegative increasing function with $f(0)=0$ such that $\n(u)=f^{-1}(u)$ is convex and $0\le f(u)\le c_p u^p$ for some $0<p\leq 1$\,. In the next statement, when $\eta=1$ then $\left\|u\right\|_{C^{\eta}}$ denotes the Lipschitz norm of $u$.

\begin{thm}[H\"older continuity up to the boundary]\label{thm.regularity.2}
Under  assumptions of Theorem \ref{thm.regularity.1}(ii),
assume in addition that $\AI$ satisfies (K2)
and that $\n=f^{-1}$ is positive convex function with $F(0)=0$. Suppose that $f\in C^\beta(\RR)$ for some $\beta>0$, that
$0\le f(a)\le c_p a^p$ when $0 \leq a \leq 1$ for some $0<p\leq 1$,
and that for any $x \in \Omega$ it holds $\|D\|_{C^\beta(B_{3r_x/2}(x))}\leq r_x^{-(2s+\beta)}$, where $r_x:=\dist(x,\partial\Omega)/2$. Then $u$ is H\"older continuous up to the boundary. More precisely, there exists  a constant $C>0$ such that
\begin{equation}\label{thm.regularity.2.bounds}
\left\|u\right\|_{C^{\eta}(\overline{\Omega})}\le C\qquad \forall\,\eta \in (0,\gamma]\cap(0,2s).
\end{equation}
\end{thm}
\noindent\textbf{Remark. } Since $u(t,x)\lesssim \Phi_1(x)\asymp\dist(x,\partial\Omega)^{\gamma}$, we get the sharp $\gamma$-H\"older continuity whenever $\gamma<2s$
(see also the remark at the end of the proof).

\noindent\textit{Proof of Theorem \ref{thm.regularity.2}. }
Fix $\eta \in (0,\gamma]$ with $\eta<2s$.
Given $x_0 \in \Omega$, set $r=\textrm{dist}(x_0,\partial\Omega)/2$
and define
\[
u_r(x):=r^{-\eta }\,u\left(x_0+rx\right)\,.
\]
With this definition, we see that $u_r$ satisfies the equation
\begin{equation}\label{eq.riscalata.f}
\A_r u_r(x)=r^{2s-\eta}f(u(x_0+rx))=r^{2s-\eta}f\left(r^{\eta} u_r(x)\right)
\end{equation}
in the rescaled domain $\Omega_r:=(\Omega-x_0)/r$, where
\begin{equation}\label{op.scaled}
\begin{split}
&\A_r w(x)=P.V.\int_{\RR^N}\big(w(x)-w(y)\big)K_r(x,y)\dy+D_r(x)w(x),\\
&K_r(x,y):=r^{N+2s}K(x_0+rx,x_0+ry),\qquad D_r(x):=r^{2s}D(x_0+rx)\,.
\end{split}
\end{equation}
Since  $u\leq C \dist(\cdot,\partial\Omega)^\gamma\leq C\dist(\cdot,\partial\Omega)^\eta$ (by Theorem \ref{Thm.upper.ell} and because $\eta \leq \gamma$), $|D(x)|\leq  C\dist(\cdot,\partial\Omega)^{-2s}$, and $\|D\|_{C^\beta(B_{3r_x/2}(x))}\leq r_x^{-(2s+\beta)}$,  we have
\[
0\leq u_r(x)\leq C \qquad\mbox{for all $x\in B_1(0)$,}\qquad \|D_r\|_{C^\beta(B_{3/2})}\leq C,
\]
with a constant $C>0$ that is independent of $r$ and $x_0$.
Also, thanks to the assumption $0\le f(u)\le c_p u^p$, for $r\leq 1$ we have (recall that $\eta<2s$, thus $\eta(1-p)<2s$)
\[
r^{2s-\eta}f\left(r^{\eta} u_r(x)\right) \le c_p r^{2s-\eta(1-p)} u_r^p(x) \le C_p \qquad\mbox{for all $x\in B_1$\,.}
\]
Furthermore, using again that $u(x)\leq C \dist(x,\partial\Omega)^{\eta}$, we see that
\[
u_r(x)\leq C\bigl(1+|x|^{\eta}\bigr) \qquad \text{for all  $x\in \RR^N$}.
\]
Thus, since $\eta<2s$ by assumption, the tails
of $u_r$ will not create any problem. Indeed, for any $x \in B_1$,
\begin{equation}\label{growthVSintegral}
\int_{\mathbb R^N\setminus B_2}u_r(y)K_r(x,y)^{-(N+2s)}\,dy\leq C\int_{\mathbb R^N\setminus B_2}|y|^\eta |y|^{-(N+2s)}\,dy\leq \bar C_0,
\end{equation}
where $\bar C_0$ is independent of $r$. As a consequence, if we localize the equation as in \eqref{localization} and apply Theorem \ref{thm.regularity.1}(ii), we get
\begin{equation*}%\label{eq:interior ur}
\|u_r\|_{C^{2s+\beta'}(B_{1/2})} \leq C\,.
 \end{equation*}
Since $\eta\leq 2s+\beta'$ (recall that $\eta<2s$), it follows that
$$
\|u_r\|_{C^{\eta}(B_{1/2})} \leq \|u_r\|_{C^{2s+\beta'}(B_{1/2})} \leq C.
$$
Noticing that
$$
[u_r]_{C^{\eta}(B_{1/2})}
=
[u]_{C^{\eta}(B_{r}(x_0))},
$$
and that $x_0$ is arbitrary, we deduce the desired H\"older continuity up to the boundary.
\qed

%%%%%%%%%%%%%%%%%%%%%%%%%%%%%%%%%%%%%%%%%%%%%%%%%%%%%%%%%%%

\noindent\textbf{Remark. } In the regime $\gamma\geq 2s$, one can actually improve the H\"older exponent through an iteration procedure. More precisely, once we know that $u\in C^\eta(\Omega)$, we can fix $|h|=1$ and consider the function
$$
v_r^h(x):=\frac{u(rx+rh)-u(rx)}{r^\zeta},
$$
where $\zeta\in (\eta,\gamma)$ has to be chosen.
Thanks to the bound $u\lesssim \dist(\cdot,\partial\Omega)^\gamma$
and the $\eta$-H\"older continuity of $u$,
we can bound $v_r^h(x)$ for $|x|\ge 1$ as
$$
|v_r^h(x)|\leq Cr^{-\zeta}\min\{|rx|^{\gamma},r^\eta\}.
$$
In particular, for any fixed $|x|\ge 1$,
$$
r^{-\zeta}\min\{|rx|^{\gamma},r^\eta\} \leq \sup_{t>0}\min\{|x|^{\gamma}t^{\gamma-\zeta},t^{\eta-\zeta}\}=|x|^{\frac{\zeta-\eta}{\gamma-\eta}\gamma},
$$
which proves that
$$
|v_r^h(x)|\leq C\left(1+|x|^{\frac{\zeta-\eta}{\gamma-\eta}\gamma}\right)\qquad \forall\,x \in \mathbb R^n.
$$
In particular, provided $\zeta>\eta$ is such that $\frac{\zeta-\eta}{\gamma-\eta}\gamma<2s$, $\zeta\leq 2s+\beta'$, and
$$
r^{2s-\zeta}|f(u(rx+rh))-f(u(rx))| \leq C \qquad \forall\, x \in B_1,
$$
(this is the case provided $2s-\zeta+\gamma p\geq 0$),
we can repeat the argument in the proof of Theorem \eqref{thm.regularity.2} to deduce that $u \in C^\zeta(\Omega)$.
In particular, as long as $\zeta\leq \min\{2s+\gamma p,2s+\beta'\}$, we can iterate this argument and actually reach any exponent strictly less than $\gamma$.
Since $2s+\beta'$ can be any exponent strictly less than $2s+2s\frac{\beta}{1-\beta}=\frac{2s}{1-\beta}$,
one deduces that $u \in C^\zeta$ with
$\zeta<\min\{2s+\gamma p,\frac{2s}{1-\beta},\gamma\}$ (in particular, when $f(a)=a^p$ then $\beta=p$ and $\zeta<\min\{2s+\gamma p,\frac{2s}{1-p},\gamma\}$).
We leave the details to the interested reader.

%%%%%%%%%%%%%%%%%%%%%%%%%%%%%%%%%%%%%%%%%%%%%%%%%%%%%%%%%%%%%%%%%%%%%%%%%%%%%%%%%%%%%%%%%%%%%%%%%%%%
%%%%%%%%%%%%%%%%%%%%%%%%%%%%%%%%%%%%%%%%%%%%%%%%%%%%%%%%%%%%%%%%%%%%%%%%%%%
\bigskip

\noindent{\sc Acknowledgments. } M.B. and J.L.V. are partially funded by MTM2014-52240-P (Spain). A.F. has been supported by ERC Grant ``Regularity and Stability of Partial Differential Equations (RSPDE)''. M.B. and J.L.V. would like to acknowledge the hospitality of the Mathematics Department of the University of Texas at Austin, where part of this work has been done. Also, M.B. is grateful for the hospitality of FIM at ETH Z\"urich. We would like to thank the anonymous referee for useful comments.

\end{document}